\documentclass{workingpaperJJB} 
\usepackage[normalem]{ulem}
\usepackage{natbib}
\usepackage{amsmath}
 \bibpunct[, ]{(}{)}{,}{a}{}{,}%
 \def\bibfont{\small}%

 \TheoremsNumberedThrough     
 \ECRepeatTheorems

 \EquationsNumberedThrough    

\usepackage{amsfonts}
\usepackage{amsmath}
\usepackage{amssymb}
\usepackage{changepage}
\usepackage{graphicx}
\usepackage{multirow}
\usepackage{enumerate}
\usepackage{booktabs}
\usepackage{color}
\usepackage{filecontents}
\usepackage{caption}
\usepackage{subcaption}
\usepackage{threeparttable}
\usepackage{tikz}
\usepackage{bm}
\usepackage[]{algorithm2e}
\RestyleAlgo{ruled}

    \newcommand{\fignote}[2]{\begin{center}\parbox[c]{#1}{\footnotesize  #2} \end{center}}
    \newcommand{\tabnote}[2]{\begin{center} \parbox[c]{#1}{\footnotesize #2} \end{center}}

\usepackage{endnotes}

\def\sym#1{\ifmmode^{#1}\else\(^{#1}\)\fi}

\usepackage{endnotes}

\DeclareMathOperator*{\minimize}{minimize}

\newcommand{\bA}{ \mathbf{A} }

\newcommand{\bb}{ \mathbf{b} }
\newcommand{\bC}{ \mathbf{C} }

\newcommand{\bw}{ \mathbf{w} }

\newcommand{\by}{ \mathbf{y} }

\newcommand{\bzero}{ \mathbf{0} }

\newcommand{\balpha}{ \boldsymbol{\alpha} }

\begin{document}


\RUNAUTHOR{Smith and Boutilier}

\RUNTITLE{Gap-Gradient Methods for Inverse Optimization}

\EquationsNumberedThrough    

\TITLE{Gap-gradient methods for solving generalized mixed integer inverse optimization: an application to political gerrymandering}

\ARTICLEAUTHORS{
\AUTHOR{Ari Smith}\AFF{Department of Industrial and Systems Engineering, University of Wisconsin-Madison, \EMAIL{ajsmith44@wisc.edu}}
\AUTHOR{Justin J. Boutilier}\AFF{Telfer School of Management, University of Ottawa, \EMAIL{boutilier@telfer.uottawa.ca}}}
\HISTORY{\today}

\ABSTRACT{

Inverse optimization has received much attention in recent years, but little literature exists for solving generalized mixed integer inverse optimization. We propose a new approach for solving generalized mixed-integer inverse optimization problems based on sub-gradient methods. We characterize when a generalized inverse optimization problem can be solved using sub-gradient methods and we prove that modifications to classic sub-gradient algorithms can return exact solutions in finite time. Our best implementation improves solution time by up to 90\% compared to the best
performing method from the literature. We then develop custom heuristic methods for graph-based inverse problems using a combination of graph coarsening and ensemble methods. Our heuristics are able to further reduce solution time by up to 52\%, while still producing near-optimal solutions. Finally, we propose a new application domain – quantitatively identifying gerrymandering – for generalized inverse integer optimization. We apply our overall solution approach to analyze the congressional districts of the State of Iowa using real-world data. We find that the accepted districting marginally improves population imbalance at the cost of a significant increase in partisan efficiency gap. We argue that our approach can produce a more nuanced data-driven argument that a proposed districting should be considered gerrymandered.

}

\HISTORY{\today}

\maketitle

\vspace{-2em}
\section{Introduction} \label{intro}

Inverse optimization (IO) has received significant attention in recent years; see \cite{chan2023inverse} for a thorough review. At a high level, inverse optimization infers unknown parameters (e.g., cost vector, constraint matrix, etc.) of a \emph{forward optimization problem} (FOP) such that an observed solution is rendered optimal (or nearly optimal) for the FOP. Inverse optimization has been applied to a wide range of problems including radiation therapy \citep{chan2014generalized}, transportation \citep{patriksson2015traffic}, and cellular biology \citep{zhao2015learning}.

The IO literature can be partitioned into two general approaches, depending on the problem setting. The first approach (``classical IO") is to assume that the observed solution is optimal for the FOP and then leverage results from duality theory to design an appropriate IO problem \citep{ahuja2001inverse, schaefer2009inverse}. The second approach (``generalized IO") is to make no assumptions about the optimality of the observed solution (i.e., it may not be optimal or even feasible) and then use a loss function to design an appropriate IO problem \citep{chan2014generalized}. 

To date, most of the IO literature has focused on convex FOPs with relatively little research on purely integer or mixed-integer FOPs. In the context of integer IO problems, there are two streams of literature. First, \cite{schaefer2009inverse} and \cite{lamperski2015polyhedral} characterize the polyhedral representation of purely integer and mixed-integer IO problems, respectively. Second, \cite{wang2009cutting} propose a cutting plane approach for efficiently solving mixed-integer IO problems. Recent research has extended this cutting plane approach to be more efficient using parallel computing \citep{wang2013branch} or trust regions \citep{bodur2022inverse}, and to solve generalized mixed-integer IO problems \citep{moghaddass2020inverse}. 

In this paper, we propose a new approach for solving generalized mixed-integer IO problems based on sub-gradient methods. Our method leverages the observation that the generalized IO problem is equivalent to minimizing a loss function over a bounded domain and many commonly used loss functions permit sub-gradient calculation. We present the conditions for when a generalized IO problem can be solved using sub-gradient methods and we test various implementations, including gradient descent and the Frank-Wolfe method. 

\subsection{Motivating example: political gerrymandering}


In the United States of America, representation in federal and state legislatures is  apportioned by democratic elections in geographic sub-regions known as legislative districts. Every 10 years, a nationwide census is conducted and the data gathered is used to redraw new legislative districts. This redrawing process, commonly referred to as \emph{political districting}, is regulated by federal and state laws. The United States Constitution requires that political districts must be of equal population, and state constitutions enforce additional requirements upon district design. For example, the state of Iowa requires that congressional districts do not split up counties (Iowa Constitution Article III \S 37) and that districts must be as compact as practicable (Iowa State Code, Chapter 42.4.4).


Political districting is often politically contentious, as there are frequent allegations -- in high-profile court cases and in the broader public political sphere -- that districts are drawn to create imbalanced voting power across populations \citep{AP2022, APNY2022, AP2019}. When such imbalance is used to privilege the representation of one political party over another, it is referred to as \emph{partisan gerrymandering}. Allegations of partisan gerrymandering often result in trials, where judges are tasked with deciding if a particular political districting is illegal, and as such must be redrawn. In these cases, the judicial system has indicated a desire for quantitative methods of measuring partisan gerrymandering, allowing for judicial rulings to be made with more confidence (and based on defensible data) (Vieth v. Jubelirer 2004). However, quantitatively identifying gerrymandering is challenging because public and judicial discourse around gerrymandering typically invokes common sense notions --- particularly of \textit{fairness} and \textit{compactness} --- which are often not quantified or critically elaborated. 

The foundation of quantitative methods for identifying political gerrymandering relies on the development of district-level metrics that quantify and align with sociopolitical values; either those legally required (e.g., compactness in many states) or those in the public discourse (e.g., ``fair" representation). Two constitutional requirements -- contiguity and equal population -- are easily quantified. However, a frequent third requirement -- compactness -- has been subject to much debate. For example, \cite{young1988measuring} demonstrated that the common sense notions of compactness should not be expected to match a single quantitative metric, because nearly (if not) all compactness metrics can be met with counterexamples that are quantitatively compact but do not feel compact to qualitative perceptions, or vice versa. Similarly, ``fair" partisan representation has been quantified with a multitude of metrics, under such names as partisan bias \citep{grofman2007future}, competitiveness of elections \citep{nagle2017competitive}, and the efficiency gap of wasted votes \citep{stephanopoulos2015partisan}.

Initial research on identifying political gerrymandering invoked hard cutoffs for a single (chosen) metric as an indicator of gerrymandering \citep{stephanopoulos2015partisan}. 
However, a univariate approach fails to consider the contingency of what values are actually attainable for real political geographies of state populations, or that a multitude of competing metrics may inhibit each other from simultaneously meeting some benchmark, even if they all align with different democratic values that are important to courts. To address this shortcoming, researchers have used computational methods to develop Markov Chain Monte Carlo approaches that create a set of ``reasonable districtings" and corresponding probability distributions for each districting metric. A proposed political districting is then compared against this information to determine whether or not it is a statistical outlier; if so, this may indicate that a districting prioritizes an undemocratic value and should be considered gerrymandering \citep{duchin2019, deford2020}.

In this paper, we propose an inverse optimization approach to quantitatively identify partisan gerrymandering. By formulating the process of drawing political districtings as a multi-objective optimization problem, inverse optimization allows us to consider how a proposed districting enacts a prioritization of some values over others with respect to what is possible in a given political geography. Thus, even if metrics that are motivated to track some notion of partisan unfairness are deemed fictitious or not a compelling legal argument in themselves, an inverse optimization analysis may show that a districting prioritizes unfairness specifically to the detriment of other constitutionally mandated objectives, such as creating compact districts or reducing population imbalance. We argue that this approach can produce a more nuanced data-driven argument that a proposed districting should be considered gerrymandered.

\subsection{Contributions}
We summarize our contributions as follows:
\begin{enumerate}
    \item We propose a new approach -- \emph{gap-gradient methods} -- for generalized inverse integer optimization that leverages the problem's similarity to minimizing a bounded convex function where querying the function gradient is computationally difficult. We characterize when a generalized IO problem can be solved using sub-gradient methods and we prove that modifications to classic sub-gradient algorithms can return exact solutions in finite time.
    We evaluate our methods using a set of instances from the MIPLIB 2017 library and we find that our best implementation is able to improve solution time by up to 90\%, compared to the best performing method from the literature.
    \item We develop a custom heuristic method for graph-based inverse problems using a combination of coarsening methods from the field of graph partitioning and ensemble methods from the field of machine learning. Our methods use ensembles of smaller problem instances to produce approximate solutions that provably converge towards the optimal inverse solution as ensemble size increases, while saving computational expense and allowing for tractable run times for large real-world problems. Our heuristic is able to reduce the median solution time by 52\%, while still producing near-optimal solutions.
    \item We propose a new application domain -- quantitatively identifying gerrymandering -- for generalized inverse integer optimization. We apply our overall solution approach to analyze the congressional districts of the State of Iowa using real-world data. We show that accepted and previously rejected district plans both greatly prioritize minimizing population imbalance over district compactness or partisan efficiency gap. However, the increased priority on population imbalance in the accepted district plan results in a 16\% decrease in population imbalance while creating a 372\% increase in partisan efficiency gap compared to the rejected plan.

\end{enumerate}
Many of the graph-based IO methods formulated and implemented in this paper (Contribution 2) can also be used in several other applications, such as flexible job shop scheduling, employee scheduling, and facility location, particularly in settings where there are competing objectives to be considered (e.g., distribution of polling locations).

\section{Literature Review}\label{litreview}

Our work contributes to three primary streams of literature: inverse optimization (Section~\ref{lit:IO}), optimization approaches for political districting (Section~\ref{lit:OptPD}), and quantitative methods for identifying gerrymandering (Section~\ref{lit:Ger}). 

\subsection{Inverse Optimization}\label{lit:IO}

The idea of IO was proposed by \cite{ahuja2001inverse} who studied linear FOPs in a classical IO setting. Since then, there has been much research on IO, including for conic FOPs \citep{iyengar2005inverse}, general convex FOPs \citep{zhang2010inverse}, purely integer FOPs \citep{schaefer2009inverse}, and mixed-integer FOPs \citep{wang2009cutting, bulut2021complexity}. See \cite{chan2023inverse} for a recent review.

In recent years, there have been several extensions to IO, often inspired by machine learning and data-driven optimization. Importantly, \cite{chan2014generalized} developed the generalized IO paradigm, where no assumptions are made about the optimiality or feasibility of the observed solutions. Other extensions include IO with multiple observations \citep{keshavarz2011imputing,babier2021ensemble}, IO with uncertain data \citep{ghobadi2018robust, aswani2018inverse}, and goodness of fit measures for IO \citep{chan2019inverse}. 

Our work is most closely related to the IO literature where the FOP is a mixed-integer linear program. In this vein, \cite{wang2009cutting} proposed the first algorithmic approach (a cutting plane algorithm) for solving mixed integer IO problems. Since then, there have been several extensions that seek to improve tractability  \citep{bodur2022inverse,moghaddass2020inverse,wang2013branch}. Most similar to our paper is the work of \cite{moghaddass2020inverse} who develop the first approach for generalized inverse mixed integer optimization by leveraging \cite{wang2009cutting}'s cutting plane method, while assuming a loss function where the optimality gap is formulated as an absolute gap. Recently (and in parallel to this work), \cite{scroccaro2023learning}  develop descent methods for learning solutions to mixed integer IO problems. However, they use a more general augmented suboptimality loss function that does not necessarily permit provable finite time convergence to the optimal solution.

We contribute to this literature in three key ways: (i) we propose a new solution method for solving generalized inverse multi-objective mixed integer optimization in finite time that uses first order methods of minimizing a bounded convex function, (ii) we propose a set of heuristic approaches inspired by ensemble methods in machine learning that are tailored to the setting where the FOP is graph-based (i.e., amenable to graph coarsening), and (iii) we propose a new and potentially impactful application domain for generalized inverse integer optimization.

\subsection{Optimization approaches to political districting} \label{lit:OptPD}

Political districting has a long history as an optimization problem in the operations research community. See \cite{ricca2013political} or \cite{swamy2022multiobjective} for a recent literature review. The literature in this area typically uses mixed-integer programming to design political districtings \citep{hess1965nonpartisan}. Most approaches focus on nonpolitical metrics such as population imbalance, contiguity, and compactness \citep{garfinkel1970optimal}. Since these problems are very computationally challenging, research has focused on developing efficient formulations and both exact and heuristic solution methods \citep{mehrotra1998optimization,validi2022imposing}. 

More recently, researchers have developed districting models that include or optimize for notions of political fairness such as efficiency gap, partisan symmetry, and competitiveness \citep{swamy2022multiobjective, king2015efficient}. For example, \cite{swamy2022multiobjective} produce a formulation that considers districting design as a multi-objective optimization problem between district compactness and variety of fairness-related metrics, and uses the formulation to explore the Pareto frontier of tradeoffs between pairs of metrics over a state's political geography. While we do not directly contribute another districting formulation, we make minor modifications to \cite{swamy2022multiobjective}'s formulation and their multi-objective approach, and use it as the FOP for our generalized inverse optimization model, demonstrating another use of their formulation.

\subsection{Quantitative Methods for Identifying Gerrymandering} \label{lit:Ger}

Initial research on identifying political gerrymandering invoked hard cutoffs for a single (chosen) metric as an indicator of gerrymandering; see \cite{stephanopoulos2018measure} for a review. In recent years, research has focused on how metrics can quantitatively identify gerrymandering using approaches that are not cutoff rules. For example, \cite{duchin2018gerrymandering} propose the use of Markov Chain Monte Carlo methods to produce a set of reasonable districtings.  
This set of districtings is then used to create statistical distributions of particular partisan fairness (or compactness) metrics. If the districting under scrutiny presents a partisan fairness metric that is a strong statistical outlier compared to the distribution, then an argument can be made that the districting explicitly prioritizes disproportionate representation (i.e., partison ``unfairness"). Note that this line of argument requires the assumption that the constructed set of districtings accurately represents average or permissible districtings. These methods have been applied to the study of various underlying metrics for district design, including partisan seat share \citep{duchin2019}, partisan symmetry \citep{deford2020}, and proportion of majority-minority districts \citep{becker2021}.

In contrast to previous approaches, we contribute a new method for quantitatively measuring how district design impacts partisan fairness. 
Rather than comparing a districting to a constructed distribution of permissible districtings like the Markov Chain Monte Carlo methods, our inverse optimization methods analyze how a district compares to the boundary of what is feasible, while critically recognizing that multiple metrics can be in competition with each other. For example, if a districting that meets a certain benchmark for one metric limits the best possible performance on other metrics, then producing judgements based on the evaluation of a single metric may be unhelpful when multiple metrics are valued. With inverse optimization, we can understand how a districting choice enacts a prioritization of metrics and their corresponding democratic values.  For example, if one argued that high political unfairness is not sufficient to disqualify a districting, then one may show that the districting under scrutiny prioritizes political unfairness explicitly to the detriment of another legal right, such as district compactness or population balance. We beleive our approach may provide a more robust argument that a districting should be disqualified.

\section{Models} \label{s.models}

In this section, we introduce the general structure of the forward multi-objective mixed-integer linear optimization problem (Section~\ref{s.fop}) and the corresponding inverse optimization problem (Section~\ref{s.inv})

\subsection{The Forward Problem}\label{s.fop}

To apply inverse optimization, we must first define the structure of the forward optimization problem (FOP). Inverse optimization implicitly assumes that the observed solutions were generated by a decision making process similar to the FOP. In the context of political districting, the decision making process is typically represented by a multi-objective mixed-integer linear optimization problem \citep{swamy2022multiobjective}.

We first introduce the general problem structure and provide the specific FOP for political districting in Section~\ref{s.res}. Let $\by\in\mathbb{R}^n\times\mathbb{Z}^{n-q}$, $q={0,1,2,...,n}$ represent the decision variables for the FOP, and let $\bA\in\mathbb{R}^{m\times n}$ and $\bb\in\mathbb{R}^m$ represent the constraints. We let $\mathcal{K}$ denote the set of objective functions. Then, the general multi-objective mixed-integer linear optimization problem can be written as
\allowdisplaybreaks
\begin{subequations} \label{generalFP}
\begin{align}
    \minimize_{\by} \quad &\balpha^T\bC \by\\
    \text{subject to} \quad &A \by \leq \bb, \\
    &\by\in\mathbb{R}^n\times\mathbb{Z}^{n-q},
\end{align}
\end{subequations}
where the rows of $\bC\in\mathbb{R}^{|\mathcal{K}|\times n}$ represent different linear objective functions and $\balpha\in\mathbb{R}^{|\mathcal{K}|}$ denotes the cost vector. Let $\mathcal{S}=\{A \by \leq \bb, \by\in\mathbb{R}^n\times\mathbb{Z}^{n-q}\}$ denote the feasible region of the FOP, let $\mathcal{B}(\mathcal{S}) = \{\bC\by |\; \by\in \mathcal{S}\}$ denote the  set of feasible sub-objective values (also known
as the FOP objective feasible space), and let $\mathcal{F}(\balpha,\mathcal{S})$ represent the set of optimal solutions. Without loss of generality, we assume that $\mathcal{S}$ is non-empty \citep{chan2014generalized}.

\subsection{Inverse Optimization Model}\label{s.inv}

In its classical form, inverse optimization seeks to, given a feasible region for a forward problem ($\mathcal{S}$) and an observed solution ($\hat{\by}$), return an objective function to the forward problem ($\bar{\balpha}$) for which the given solution is optimal -- a property called \emph{inverse feasibility} \citep{chan2023inverse}. Given a feasible solution to the forward problem, $\hat{\by}\in\mathcal{S}$, we let $\mathcal{C}(\hat{\by},\mathcal{S}) = \{\balpha\in\mathbb{R}^{|\mathcal{K}|}\; | \; \hat{\by}\in\mathcal{F}(\balpha,\mathcal{S})\}$ denote the inverse-feasible region, which includes all cost vectors $\balpha$ that render $\hat{\by}$ optimal.

In data-driven contexts, it may not be reasonable to assume that $\hat{\by}$ is optimal for any possible objective (i.e., $\mathcal{C}(\hat{\by},\mathcal{S})= \varnothing$) or even feasible for the FOP (i.e., $\hat{\by}\not\in\mathcal{S}$), so instead we seek to minimize some loss function of the forward objective that penalizes the extent to which inverse-feasibility is not satisfied. 
We denote the loss function $\ell(\hat{\by}, \mathcal{S}, \balpha)$, and note that if a solution $\hat{\by}$ does have some classically inverse feasible objective, then for any $\balpha \in \mathcal{C}(\hat{\by}, \mathcal{S})$, we have $\ell(\hat{\by}, \mathcal{S}, \balpha) = 0$. \cite{chan2023inverse} enumerate five possible loss functions that are useful in the data-driven inverse optimization context.

For any such loss function, the corresponding inverse optimization problem seeks to minimize the loss function, with constraints placed on the allowable objectives $\balpha$, i.e., $\balpha\in\mathcal{A}$. We can write this problem as:
\begin{align*}
    \minimize_{\bm{\alpha}} \quad & \ell(\hat{\by}, \mathcal{S}, \balpha)  \tag{GIO$(\hat{\by}, \mathcal{S})$}\\
    \text{subject to } \quad & \balpha \in \mathcal{A}.
\end{align*}

Many choices of loss functions have inverse optimization formulations that can be more practically expressed. For example, for the absolute sub-optimality loss function $\ell_{\text{ABS}}(\hat{\by}, \mathcal{S}, \balpha) = \balpha^\intercal \bC \hat{\by} - \underset{\by \in \mathcal{S}} {\min}\:\balpha^\intercal \bC \by$, which measures the difference in FOP objective value between the input $\hat{\by}$ and an optimal FOP solution for the given value of $\balpha$, the equivalent inverse optimization formulation can be written as:
\begin{align*}
    \minimize_{\bm{\alpha},\; \xi_{\text{ABS}}} \quad & \xi_{\text{ABS}}  \tag{GIO$_{\text{ABS}}(\hat{\by}, \mathcal{S})$}\\
    \text{subject to } \quad & \bm{\alpha}^\intercal \bC \hat{\by} \leq \bm{\alpha}^\intercal \bC \by + \xi_{\text{ABS}}, \quad \forall \: \by \in \mathcal{S},\\
    & \balpha \in \mathcal{A}. 
\end{align*}

One potentially useful choice of $\mathcal{A}$ for many applications is the unit simplex in dimension $|\mathcal{K}|$, because then $\balpha$ may be interpreted as a convex combination of weights for multiple objectives, allocating fractional amounts of importance to each objective in a form that adds up to 1. We also note that when $C$ is the identity matrix, the generalized multi-objective inverse optimization problem GIO$(\hat{\by}, \mathcal{S})$ becomes a standard generalized inverse optimization problem, and all following solution approaches remain applicable without loss of generality.

\section{Solution Approaches}

In this section, we recap existing cutting plane methods (Section~\ref{sec:cutplane}) and introduce a new approach for solving inverse mixed integer optimization problems (Section~\ref{sec:gap-gradient}).

\subsection{Cutting Planes}\label{sec:cutplane}

The first cutting plane algorithm for inverse mixed integer linear optimization was tailored to the classical inverse case that assumes the existence of an inverse-feasible objective \citep{wang2009cutting}. The algorithm alternates between a master problem and a separation problem. The master problem produces an optimal forward problem objective ($\balpha$) that is hypothesized to be inverse feasible (i.e., $\balpha\in \mathcal{C}(\hat{\by},\mathcal{S})$) under partial knowledge of the forward feasible region ($\mathcal{S}$). The separation problem solves the forward optimization problem with the hypothesized objective ($\balpha$) to produce an extreme point that can be added to the set $\mathcal{S}$. With each new extreme point that is added to $\mathcal{S}$, the set of potentially inverse-feasible objectives $\mathcal{C}(\hat{\by},\mathcal{S})$ decreases in size. When the forward problem returns an extreme point that is already known to be in $\mathcal{S}$, further iterations will not remove elements from $\mathcal{C}(\hat{\by},\mathcal{S})$, and so the output of the master problem is known to be an inverse-feasible solution. Thus, the algorithm returns the solution to the master problem and terminates. Since the forward problem is a mixed integer linear program, it is guaranteed to have a finite number of extreme points, and so the algorithm is guaranteed to (eventually) terminate.

More recently, \cite{moghaddass2020inverse} proposed an extension to Wang's algorithm that is suitable for data-driven generalized inverse optimization, in the case of an absolute sub-optimality loss function. As such, they use GIO$_{\text{ABS}}(\hat{\by}, \mathcal{S})$ as their master problem, and alternate between the master problem and the separation problem until an objective is found that minimizes the loss function. We present their algorithm structure (Algorithm~\ref{moghaddass}) using our established notation in \ref{sec:EC_algstruc}. In \ref{sec:EC_MTrelgap}, we contribute a minor modification to the algorithm proposed by \cite{moghaddass2020inverse} that utilizes the relative sub-optimality loss function, rather than absolute sub-optimality.

\subsection{Gap-Gradient Methods}\label{sec:gap-gradient}

In this section, we describe our proposed solution methods -- \emph{gap-gradient methods} -- and highlight various extensions. Our methods are applicable to a wide range of loss functions and constraint sets for the cost vector.

We first define what we call the \textit{gap function} for a given generalized inverse optimization problem. The gap function is simply the loss function values over the feasible region of $\balpha$, which we formally defined as follows.

\begin{definition}[Gap Function]
     For a generalized inverse optimization problem with fixed $\hat{\by}$ and $\mathcal{S}$, and some loss function $\ell(\hat{\by}, \mathcal{S}, \balpha)$, the gap function is defined as $\xi(\balpha) = \ell(\hat{\by}, \mathcal{S}, \balpha)$, $\forall \balpha\in\mathcal{A}$. 
\end{definition}

Note that the gap function will be a convex function as long as the loss function is a convex function of $\balpha$, given fixed $\hat{\by}$ and $\mathcal{S}$. For the example where the loss function is absolute sub-optimality and $\balpha$ is constrained to the unit simplex, the gap function is $\xi_{ABS}(\balpha) = \balpha^\intercal \bC \hat{\by} - \underset{\by \in \mathcal{S}} {\min}\:\balpha^\intercal \bC \by$, where $\balpha\in\{||\bm{\alpha}||_1 = 1, \bm{\alpha} \geq 0\}$. To keep the rest of the paper succinct, all discussion henceforth concerns the case in which the loss function $\ell$ is absolute sub-optimality and the gap function domain $\mathcal{A}$ is the unit simplex.


\begin{proposition}\label{thm:gap_abs}
    The absolute gap function $\xi_{ABS}(\balpha)$ is a convex function.
\end{proposition}

\proof{Proof of Proposition~\ref{thm:gap_abs}:}
First, note that the epigraph of $\xi_{ABS}(\balpha)$ is $\{\xi_{ABS}\in\mathbb{R} \; | \; \xi_{ABS}\geq \balpha^\intercal \bC \hat{\by} - \underset{\by \in \mathcal{S}} {\min}\:\balpha^\intercal \bC \by\}$, which is equivalent to $\{\xi_{ABS}\in\mathbb{R} \; | \; \xi_{ABS}(\balpha)\geq \balpha^\intercal \bC \hat{\by} - \balpha^\intercal \bC \by, \forall \by \in \mathcal{S}\}$. This space is the intersection of half-planes, and as such, is convex. Therefore, the gap function is a convex function because its epigraph is convex.
\halmos
\endproof

\begin{figure}
    \centering
    \begin{tikzpicture}[scale=0.9]
    \draw [->] (0, 0) -- (0, 5);
    \draw [->] (0, 0) -- (5, 0);
    \draw [pink] [fill=pink] (1, 5) -- (1, 3.4) -- (1.8, 1.7) -- (3.6, 1) -- (5, 1) -- (5,5);
    \draw [red] (1, 5) -- (1, 3.4) -- (1.8, 1.7) -- (3.6, 1) -- (5, 1);
    \draw [cyan] [->] (2.7, 1.8) -- (3.06, 2.52) node [right] {$\balpha^*$};
    \filldraw [blue] (2.7, 1.8) circle (1pt) node [right] {$\hat{\by}$};
    \node at (5, -0.4) {$c_1{}^\intercal \by$};
    \node at (-0.6, 5) {$c_2{}^\intercal \by$};
    \node at (2.5, -1) {(a)};
    \end{tikzpicture}
    \hspace{5pt}
    \begin{tikzpicture}[scale=0.9]
    \draw (0, 0) -- (5, 0);
    \draw [dashed] (0, 0) -- (0, 5);
    \draw [dashed] (5, 0) -- (5, 5);
    \draw [red] (0, 1.7) -- (1.6, 0.856917) -- (3.6, 0.4192) -- (5, 0.7);
    \filldraw [red] (0, 1.7) circle (1pt) node [above right] {$\xi_{\text{ABS}}(\balpha)$};
    \filldraw [red] (5, 0.7) circle (1pt);
    \filldraw [cyan] (3.6, 0) circle (1pt) node [below] {$\balpha^*$};
    \node at (0, 5.4) {$\xi$};
    \node at (0.1, -0.6) {\small$\alpha = \begin{pmatrix} 1 \\ 0 \end{pmatrix}$};
    \node at (5.1, -0.6) {\small$\alpha = \begin{pmatrix} 0 \\ 1 \end{pmatrix}$};
    \node at (2.5, -1) {(b)};
    \end{tikzpicture}
    \hspace{0pt}
    \caption{An example inverse optimization problem with two sub-objectives. (a) The FOP objective feasible space ($conv(\mathcal{B})$),  inverse input ($\hat{\by}$), and the inverse solution ($\balpha^*$), and (b) the corresponding absolute gap function ($\xi_{\text{ABS}}$) and its minimizer ($\balpha^*$).}
    \label{fig:gap_space_example}
\end{figure}
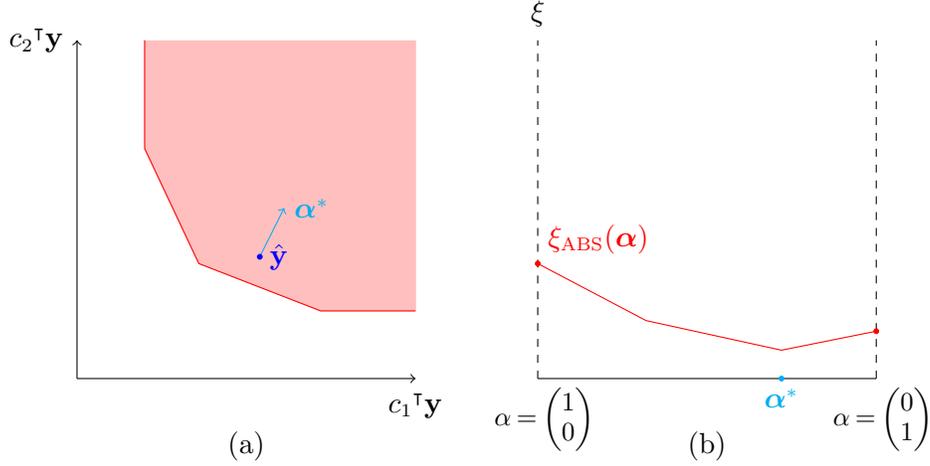

Figure~\ref{fig:gap_space_example} depicts an example inverse optimization problem with two sub-objectives, $c_1$ and $c_2$. Figure~\ref{fig:gap_space_example}(a) displays the convex hull of the set of feasible subobjective values for the forward problem ($conv(\mathcal{B})$), the inverse input ($\hat{\by}$), and the FOP multi-objective weighting ($\balpha^*$) that minimizes the absolute sub-optimality loss function $\ell_{\text{ABS}}$ given $\hat{\by}$ and $\mathcal{S}$ (i.e., the optimal inverse solution). 
Figure~\ref{fig:gap_space_example}(b) displays the corresponding absolute gap function ($\xi_{ABS}(\balpha)$). The horizontal axis of Figure~\ref{fig:gap_space_example}(b) describes the set of allowable multi-objective weightings for the two sub-objectives as a convex combination ranging from $\begin{pmatrix} 1 \\ 0 \end{pmatrix}$ to $\begin{pmatrix} 0 \\ 1 \end{pmatrix}$. Note that the vertices of $conv(\mathcal{B})$ correspond to the piecewise linear components of $\xi$, and that $\balpha^*$ (the solution to GIO$_{\text{ABS}}(\hat{\by}, \mathcal{S})$) is also the minimizer of $\xi_{ABS}(\balpha)$.

Now, consider the cutting plane algorithm used to solve GIO$_{\text{abs}}$ (e.g., \cite{moghaddass2020inverse}). At a given iteration $k$, we use $\balpha^{(k)}$ to solve the FOP and obtain a solution $\by^{(k)}$ that is optimal for objective weighting $\balpha^{(k)}$. Recall that $\mathcal{C}(\by^{(k)},\mathcal{S})$ denotes the set of $\balpha$ values that render $\by^{(k)}$ optimal for the FOP.

\begin{proposition}\label{thm:gap_func}
    If $\balpha^{(k)}\in\mathcal{C}(\by^{(k)},\mathcal{S})$, then $\xi_{ABS}(\balpha^{(k)})=\balpha^{(k)\intercal} \bC \hat{\by} - \balpha^{(k)\intercal} \bC \by^{(k)}$.
\end{proposition}
\proof{Proof of Proposition~\ref{thm:gap_func}:}
    By the definition of the absolute sub-optimality loss function, $\xi_{\text{ABS}}(\balpha^{(k)}) = \underset{\by \in \mathcal{S}}{\max} \: \balpha^{(k)\intercal} \bC \hat{\by} - \balpha^{(k)\intercal} \bC \by$. Since $\by^{(k)} = \underset{\by \in \mathcal{S}}{\argmin} \: \balpha^{(k)\intercal} \bC \by$ by the definition of $\by^{(k)}$, it must be true that $\xi_{\text{ABS}}(\balpha^{(k)}) = \underset{\by \in \mathcal{S}}{\max} \: \balpha^{(k)\intercal} C \hat{\by} - \balpha^{(k)\intercal} C \by^{(k)}$.
\halmos
\endproof

The implication of Proposition~\ref{thm:gap_func} is that we can efficiently query subgradients of the gap function, $\xi_{\text{ABS}}(\balpha^{(k)})$. We demonstrate this via two cases. First, suppose that $\by^{(k)}$ is a vertex of conv$(\mathcal{B})$, then  $\mathcal{C}(\by^{(k)},\mathcal{S})$ is a closed convex subset of $\mathcal{A}$ with dimension $\mathcal{K} - 1$ (by the definition of a vertex). For any $\balpha$ in the interior of $\mathcal{C}(\by^{(k)},\mathcal{S})$, the \emph{gradient} of $\xi_{\text{ABS}}(\balpha)$ can be calculated as $\nabla (\bm{\alpha}^\intercal C \hat{\by} - \bm{\alpha}^\intercal C \by^{(k)}) = C (\hat{\by} - \by^{(k)})$. Second, suppose that $\by^{(k)}$ is \emph{not} a vertex of conv$(\mathcal{B})$. Then we note that $\xi_{\text{ABS}}(\balpha^{(k)}) = \balpha^{(k)\intercal} C \hat{\by} - \balpha^{(k)\intercal} C \by^{(k)}$, and $\forall \balpha \in \mathcal{A}$, $\xi_{\text{ABS}}(\balpha) \geq \balpha^\intercal C \hat{\by} - \balpha^\intercal C \by^{(k)}$, since $\underset{\by \in \mathcal{S}}{\max} \: \bm{\alpha}^\intercal C \hat{\by} - \bm{\alpha}^\intercal C \by \geq \bm{\alpha}^\intercal C \hat{\by} - \bm{\alpha}^\intercal C \by^{(k)}$. Thus, the linear expression $\balpha^\intercal C \hat{\by} - \balpha^\intercal C \by^{(k)}$ defines a \emph{subtangent plane} of $\xi_{\text{ABS}}(\balpha)$, and as such the gradient of the subtangent plane, $\nabla (\bm{\alpha}^\intercal C \hat{\by} - \bm{\alpha}^\intercal C \by^{(k)}) = C (\hat{\by} - \by^{(k)})$ constitutes a \emph{subgradient} of $\xi_{\text{ABS}}(\balpha)$.

\begin{remark}\label{rem:graddesc}
    Suppose a generalized inverse optimization problem has a gap function $\xi(\balpha)$ that is convex, a gap function domain $\mathcal{A}$ that is convex, and that the gap function is lower bounded over its domain. If it is possible to query a subgradient of $\xi(\balpha)$ for any $\balpha \in \mathcal{A}$, then any subgradient method is guaranteed to converge towards the gap function minimizer as the number of iterations approaches infinity. Since we can query subgradients of $\xi_{\text{ABS}}(\balpha)$ at any $\balpha^{(k)}$, we thus have the possibility of finding the gap function minimizer with any first-order method of minimizing convex functions over bounded domains.
\end{remark}

Additional generalized inverse optimization scenarios that satisfy the conditions stated in Remark~\ref{rem:graddesc} include a relative sub-optimality loss function where $\mathcal{A}$ is the unit simplex \citep{chan2023inverse}. See \ref{sec:relgap_adaptation} for details on how our method applies to the relative sub-optimality loss function.

In general, subgradient methods do not necessarily guarantee finite-time convergence like the cutting plane method presented by \cite{moghaddass2020inverse}, which alternates solving relaxed master
problems and generating cutting planes. For example, the descent-based methods developed by \cite{scroccaro2023learning} are not guaranteed to converge upon loss function minimizers in finite time. However, in the case of a linear piece-wise gap function, each iteration of a subgradient method is capable of producing a linear cutting plane that can be added to the master problem and, when a sufficient set of cutting planes has been found, return the exact minimizer. By occasionally solving the master problem between subgradient method steps, finite-time convergence can be guaranteed. Moreover, because subgradient methods do not need to solve the master problem at every iteration, it is possible that they may converge towards the function minimizer at a faster rate than traditional cutting plane algorithms. Subgradient methods also have the ability to query subgradients with respect to non-optimal FOP solutions for heuristic approximations and stochastic descent methods, which we explore in Section~\ref{s.res}. With these ideas in mind, we implement solution methods that make use of two major methods from the literature; projected gradient descent \citep{boyd2003subgradient} and the Frank-Wolfe method \citep{frank_algorithm_1956}.

\subsubsection{Projected gradient descent (PGD).}

To execute an iteration of projected gradient descent, given some step-size $t_k > 0$, we choose our next $\bm{\alpha}$ with $\bm{\alpha}^{(k)} \leftarrow \bm{\alpha}^{(k - 1)} - t_k (C (\hat{\by} - \by^{(k - 1)}))$, and then project $\bm{\alpha}^{(k)}$ onto its domain $\mathcal{A}$ \citep{boyd2003subgradient}. For a sufficiently small step size $t_k$, if $\bm{\alpha}^{(k)}$ is not the minimizer of $\xi_{\text{ABS}}$, then $\bm{\alpha}^{(k + 1)}$ will yield a lower optimality gap. Further, there exist multiple choices of the sequence $t_k$ such that the method will converge upon an optimal solution as $k$ approaches infinity \citep{boyd2003subgradient}. 

\begin{remark}
    Once all facet-defining tangent planes that are tangent at the gap function minimizer are known, then a single run of the master problem will yield the exactly optimal solution. We can use this knowledge to create a descent method that will be guaranteed to terminate at an optimal solution rather than infinitely converge on it.
\end{remark}

We formalize this idea as follows. Let $\xi_\text{ABS}(\balpha)$ be the absolute gap function generated by a polyhedral forward feasible region $\mathcal{S}$ and an inverse input $\hat{\by}$. Let $\balpha^*$ be a minimizer of $\xi_\text{ABS}(\balpha)$. Let $\balpha^{(k)}$ be generated by a projected subgradient descent method that uses step-size sequence $t_k$, and let $\mathcal{S}^k$ denote a set of FOP solutions $\by^{(0)} \dots \by^{(k)}$ found in the descent process. At every step $k$, a subtangent plane to $\xi_\text{ABS}(\balpha)$ at $\balpha^{(k)}$ is generated from $ \by^{(k)} \in \mathcal{S}^k \supseteq \mathcal{S}^{k-1}$. If $\balpha^*$ is a unique minimizer of $\xi_\text{ABS}(\balpha)$, then it corresponds to a facet of conv$(\mathcal{B})$ that contains $\bC \hat{\by} - \xi_{\text{ABS}}(\balpha^*)\mathbf{1}$ in its interior (proof that $\bC \hat{\by} - \xi_{\text{ABS}}(\balpha^*)\mathbf{1}$ is contained on the boundary of conv($\mathcal{B}$) follows from Theorem 1 (b) of \cite{chan2014generalized}). Similarly, each facet-defining tangent plane of $\xi_\text{ABS}(\balpha)$ corresponds to a vertex of conv$(\mathcal{B})$. Let $B^k$ denote a set of points on the boundary of conv($\mathcal{B}$) corresponding to $\mathcal{S}^k$, i.e. the set $\{\bC\by^{(0)} \dots \bC\by^{(k)}\}$. Finally, note that if $\balpha^*$ is not a unique minimizer, then the above is true for a lower dimensional face of conv$(\mathcal{B})$ rather than a facet.

\begin{proposition}\label{lemma:MPsolve}
    Suppose $\balpha^*$ is a unique minimizer of $\xi_{ABS}(\balpha)$. If the convex hull of $B^k$ contains $\bC \hat{\by} - \xi_{\text{ABS}}(\balpha^*)\mathbf{1}$ in the interior of one of its facets, then solving the master problem GIO$_{\text{ABS}}(\hat\by, \mathcal{S}^k)$ will yield $\balpha^*$. 
\end{proposition}

We provide proof of this proposition in \ref{ECproofMPsolve}. As stated in the proof, if conv$(\mathcal{B})$ is a set of dimension $|\mathcal{K}|$, then there exists at least one set of only $|\mathcal{K}|$ elements of conv$(\mathcal{B})$ that need to be discovered for this condition to be met. We henceforth refer to the case where solving the master problem yields the optimal solution as the \emph{MP solve condition}.

\begin{definition}[MP Solve Condition]
    For a sequence $\balpha^{(k)}$ generated by an inverse optimization solution method, the MP Solve Condition is achieved at an iteration $k$ when solving GIO$_\text{ABS}(\hat\by, \mathcal{S}^k)$ yields $\balpha^*$. As stated in Proposition~\ref{lemma:MPsolve}, this occurs when conv$(B^k)$ contains $\bC \hat{\by} - \xi_{\text{ABS}}(\balpha^*)\mathbf{1}$ in the interior of one of its facets.
\end{definition}

Our implementation of a terminating projected gradient descent solution algorithm is presented in Algorithm~\ref{ggproj}. Our algorithm iterates through steps of projected gradient descent, and then solves the master problem whenever a descent step neither decreases the gap function by the maximum expected amount nor returns a new FOP extreme point. This occurs when $\bm{\alpha}^{(k) \intercal} C \hat{\by} - \bm{\alpha}^{(k) \intercal} \bC \by^{(k)} > \bm{\alpha}^{(k-1) \intercal} \bC \hat{\by} - \bm{\alpha}^{(k-1) \intercal} C \by^{(k-1)} + (\balpha^{(k)} - \balpha^{(k-1)})^\intercal (C (\hat{\by} - \by^{(k-1)}))$ and $\mathcal{S}^k = \mathcal{S}^{k-1}$. One possible case when this occurs is when a sufficient number of facets of $\xi_{ABS}(\balpha)$ surrounding the minimizer have been found and the descent method is repeatedly travelling between them. In this case, the MP solve condition has been reached, and the final output of GIO$_{\text{ABS}}(\hat\by, \mathcal{S}^k)$, denoted $\balpha^{\text{final}}$ will be the gap function minimizer $\balpha^*$, which can be verified by solving the FOP with $\balpha^{\text{final}}$. If the resulting gap function value is equal to the objective value of GIO$_{\text{ABS}}(\hat\by, \mathcal{S}^k)$, then $\balpha^{\text{final}} = \balpha^*$ and the algorithm terminates. If this FOP solve yields a different gap function value than the objective value returned by GIO$_{\text{ABS}}(\hat\by, \mathcal{S}^k)$, then the MP solve condition had not yet been reached. The termination criterion may have been reached by descent steps travelling between facets of $\xi_{ABS}(\balpha)$ that are not adjacent to the minimizer $\balpha^*$, and thus the descent step size is too large. In this case, the step size is divided by two and iterations of projected gradient descent proceed until the termination criterion is reached again.

\begin{algorithm}
 \SetKwInOut{Input}{input}\SetKwInOut{Output}{output}
 \Input{$C, \: \hat{\by}$, FOP, $t$}
 \Output{$\bm{\alpha}^\text{best}, \: \xi^\text{final}$}
 $k = 0, \: \mathcal{S}^k \leftarrow \emptyset, \: \bm{\alpha}^{(k)} \leftarrow \frac{1}{\mathcal{K}} \mathbf{1}$\;
 $\by^{(k)} \leftarrow \text{FP}(\bm{\alpha}^{(k)})$\;
 $\xi_\text{ABS} \leftarrow \bm{\alpha}^{(k) \intercal} C \hat{\by} - \bm{\alpha}^{(k) \intercal} C \by^{(k)}$\;
 \While{$\bm{\alpha}^{(k) \intercal} C \hat{\by} > \bm{\alpha}^{(k) \intercal} C \by^{(k)}$}{
  $k \leftarrow k + 1$\;
  $\mathcal{S}^k \leftarrow \mathcal{S}^{k-1} \cup \by^{(k-1)}$\;
  $\bm{\alpha}^{(k)} \leftarrow \bm{\alpha}^{(k-1)} - t (C (\hat{\by} - \by^{(k-1)}))$\;
  $\bm{\alpha}^{(k)} \leftarrow \text{proj}_{\Delta^\mathcal{K}}(\bm{\alpha}^{(k)})$\;
  $\by^{(k)} \leftarrow \text{FOP}(\bm{\alpha}^{(k)})$\;
  \eIf{$\balpha^{(k)} = \balpha^{(k-1)}$\textbf{or} $\bm{\alpha}^{(k) \intercal} C \hat{\by} - \bm{\alpha}^{(k) \intercal} C \by^{(k)} > \xi_\text{ABS} + (\balpha^{(k)} - \balpha^{(k-1)})^\intercal (C (\hat{\by} - \by^{(k-1)}))$ \textbf{and} $\by^{(k)} \in \mathcal{S}^k$}{
    $\balpha^{\text{final}}, \xi^{\text{final}} \leftarrow \text{GIO}_{\text{ABS}}(\hat{\by}, \mathcal{S}^k)$\;
    $k \leftarrow k + 1$\;
    $\by^{(k)} \leftarrow \text{FOP}(\bm{\alpha}^{\text{final}})$\;
    \eIf{$\balpha^{\text{final} \intercal} C \hat\by - \balpha^{\text{final} \intercal} C \by^{(k)} = \xi^\text{final}$}{
    \textbf{stop}
    }{
    $\balpha^{(k)} \leftarrow \balpha^{\text{final}}$\;
    $\mathcal{S}^k \leftarrow \mathcal{S}^{k-1} \cup \by^{(k-1)}$\;
    $t \leftarrow \frac{t}{2}$\;
    }
  }{
   $\xi_\text{ABS} \leftarrow \bm{\alpha}^{(k) \intercal} C \hat{\by} - \bm{\alpha}^{(k) \intercal} C \by^{(k)}$\;
 }
}    
 \caption{Gap-gradient projected gradient descent method.}\label{ggproj}
\end{algorithm}

\begin{proposition}\label{prop:ggtermination}
    Algorithm~\ref{ggproj} (i) terminates in a finite number of iterations, and (ii) returns the gap function minimizer.
\end{proposition}

Proof of Proposition~\ref{prop:ggtermination} is provided in Appendix~\ref{ECproofGGterminate}.

\subsubsection{Accelerated projected gradient descent (PGD-A).}
We implement Polyak's Heavy Ball Method as part of Algorithm~\ref{ggproj} by adding a momentum term with a coefficient $\beta \in (0, 1]$ to the step formulation before projecting onto the domain \citep{polyak1964some}. We use the same termination method as in the projected gradient descent method. See Algorithm~\ref{ggheavyball} in \ref{sec:EC_algstruc} for pseudocode of the algorithm structure with the addition of a momentum term, utilizing the same termination technique as projected gradient descent. 
    
\subsubsection{Frank-Wolfe method (FW).}
The Frank-Wolfe method of optimization over a convex bounded function entails querying the loss function gradient at a given point, finding the minimizer of said gradient in the function domain, and moving in the direction towards said point with a step size that decreases at a a rate of $O(\frac{1}{k})$. In our case, the domain of the function we are minimizing ($\mathcal{A}$) is a unit simplex, so the task of finding the gradient minimizer is trivial; we find the smallest component $i$ of the gradient, and the minimizer is the vector with $1$ at component $i$ and 0 for all others. For a Lipschitz-continuous convex function (e.g., the absolute gap function, $\xi_{\text{ABS}}$), the algorithm converges on an optimal solution at a rate of $O(\frac{1}{k})$ \citep{frank_algorithm_1956}. See Algorithm~\ref{ggfw} in \ref{sec:EC_algstruc} for our implementation of the Frank-Wolfe method in the context of generalized inverse optimization, using the same termination method as projected gradient descent.

\begin{remark}\label{early_stopping}
We note that it is actually possible to converge on the optimal solution using the Frank-Wolfe method without necessarily being able to query the gradient accurately, so long as one can find the point in the gap space that minimizes said gradient. For the case where the loss function is absolute sub-optimality and $\mathcal{A}$ is the unit simplex, we only need to consider $|\mathcal{K}|$ such points in our domain, so it may be possible to not solve the forward problem to completion, as long as we do enough work to be certain which of these points would be the gradient minimizer (i.e., the optimal solution will produce a gradient with a specific component being the smallest). As soon as we know this, we can abort solving the forward problem and take the next Frank-Wolfe step. Since solving the MIP forward problem is generally the most computationally complex step of all algorithms discussed so far, this has the potential to significantly improve time-per-iteration in the Frank-Wolfe algorithm (and thus potentially in overall time). However, if our method of executing early cutoffs is not able to accurately query the gap function value and yield tangent planes, we will not be able to \textit{only} conduct partial FOP solves and then terminate the algorithm with an MP solve condition. However, one may choose to run several partial solve iterations before changing the approach and using fewer full solve iterations in a close neighborhood of the solution before terminating. Note that the gap function tangent planes generated when far from the solution are not necessarily useful for terminating the master problem, and only those planes tangent to facets that include a function minimizer need to be known.
\end{remark}

We provide a formal method for operationalizing this idea in \ref{FWpartialSolve}, but we leave the practical implementation of this approach for future work.

\section{Experimental Evaluation of Solution Methods}\label{MIPLIBExperiments}

In this section, we evaluate our proposed methods and compare them to an existing method using a set of multi-objective inverse mixed integer optimization problems derived from FOPs obtained from the MIPLIB 2017 mixed integer programming library \citep{miplib}.

\subsection{Experimental Setup}

We implement and compare our three proposed solution methods: 1) projected gradient descent (PGD), 2) accelerated projected gradient descent  (PGD-A), and 3) Frank-Wolfe method (FW). We compare our approaches with the only applicable solution method from the literature: the Moghaddass-Terekhov cutting plane method (CP) \citep{moghaddass2020inverse}. For PGD and PGD-A, the initial step size is chosen so that the Euclidean norm of the first descent step is 0.1, i.e. $||\balpha^{(1)} - \balpha^{(0)}||_2 = 0.1$. For PGD-A, we use a momentum coefficient of 0.5. We note that our approach produces a conservative estimate for the performance of these methods because a more rigorous search for optimal step sizes and momentum coefficients may substantially improve solution times.

The list of selected MIPLIB instances (i.e., FOPs) is included in \ref{MIPLIBfops}. These nine instances were chosen because (1) they were flagged as ``benchmark-suitable'' and rated as ``easy'' solution difficulty, allowing for relatively tractable run times for the inverse solution algorithms, (2) they contained at least 128 continuous variables which can serve as reasonable sub-objectives, in addition to containing integer or binary variables, (3) they were neither infeasible nor unbounded for any sampled objective weights $\balpha$ in our computational trials, and (4) they span a variety FOP structures that are suited for different application settings.

To adapt each of the nine FOPs into a multi-objective generalized inverse optimization problem, we conduct the following process. For each value of $k$ in the sequence $[4, 8, 16, 32, 64, 128]$:
\begin{enumerate}
    \item We randomly sample $k$ continuous variables from the FOP formulation to be the individual sub-objectives. Each sub-objective function in our multi-objective formulation is the value of one of these sampled variables. Thus each row in the matrix $\bC$ contains all entries equalling 0 except for one sampled variable, which has coefficient 1. We execute this sampling process 3 times per FOP.
    \item For each sampled matrix $\bC$, we sample 3 different coefficient vectors $\balpha$ from the unit simplex in $k$ dimensions.
    \item For each sampled combination of $\bC$ and $\balpha$, we solve the FOP with the objective function $\balpha^\intercal \bC \by$ to obtain our inverse input solution $\hat\by$. In total, this gives us 9 inverse inputs per FOP formulation (per $k$).
    \item We used each of the four inverse solution methods to solve the inverse optimization problem for each of the nine input solutions. We use the absolute sub-optimality loss function and the unit simplex for the set of allowable objectives. For each solve, we record the total run time and the total number of iterations.
\end{enumerate}

All computational experiments were executed on an HP EliteDesk 800 G4 TWR with an Intel(R) Core(TM) i5-8500 CPU running at 3.00GHz with 6 Cores and 6 Logical Processors and 16.0 GB of RAM. All code was written and executed in Python version 3.7.3 and all FOPs were solved using Gurobi version 9.0.2. A maximum running time of 360 seconds was used in each trial before cutting off the solution algorithms.

\subsection{Results}

\begin{figure}
    \centering
    \begin{subfigure}{0.3\textwidth}
    \includegraphics[scale=0.4,trim={0 0 0 0pt},clip]{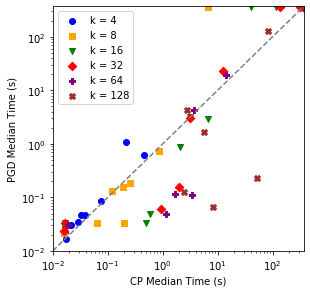}
    \caption{}\label{Exp1A}
    \end{subfigure}
    \begin{subfigure}{0.3\textwidth}
    \includegraphics[scale=0.4,trim={0 0 0 0pt},clip]{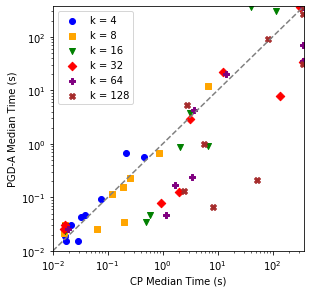}
    \caption{}\label{Exp1Ba}
    \end{subfigure}
    \begin{subfigure}{0.3\textwidth}
    \includegraphics[scale=0.4,trim={0 0 0 0pt},clip]{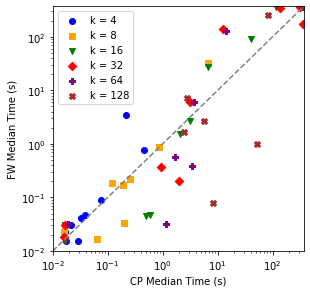}
    \caption{}\label{Exp1Bb}
    \end{subfigure}
    \caption{A comparison of median solution times for each value of $k$ (marker type): (a) PGD vs. CP, (b) PGD-A vs. CP, and (c) FW vs. CP. For each marker type (value of $k$), there are nine markers (one for each FOP instance) that denote the median solution time across all samples for that FOP instance.}
    \label{fig:MIPLIBExp1}
\end{figure}

\begin{figure}
    \centering
    \begin{subfigure}{0.3\textwidth}
    \includegraphics[scale=0.4,trim={0 0 0 0pt},clip]{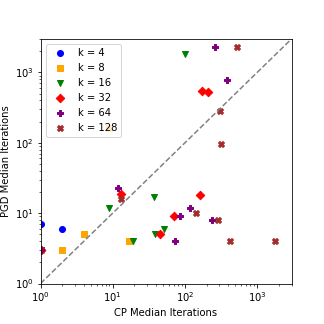}
    \caption{}\label{Exp1A2}
    \end{subfigure}
    \begin{subfigure}{0.3\textwidth}
    \includegraphics[scale=0.4,trim={0 0 0 0pt},clip]{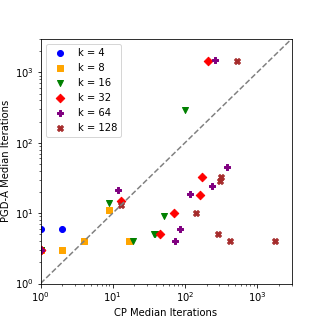}
    \caption{}\label{Exp1Ba2}
    \end{subfigure}
    \begin{subfigure}{0.3\textwidth}
    \includegraphics[scale=0.4,trim={0 0 0 0pt},clip]{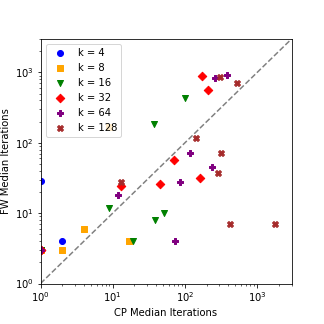}
    \caption{}\label{Exp1Bb2}
    \end{subfigure}
    \caption{A comparison of the median number of iterations for each value of $k$ (marker type): (a) PGD vs. CP, (b) PGD-A vs. CP, and (c) FW vs. CP.}
    \label{fig:MIPLIBExp2}
\end{figure}

\begin{figure}
    \centering
    \begin{subfigure}{0.48\textwidth}
        \centering
        \includegraphics[scale = 0.3, clip]{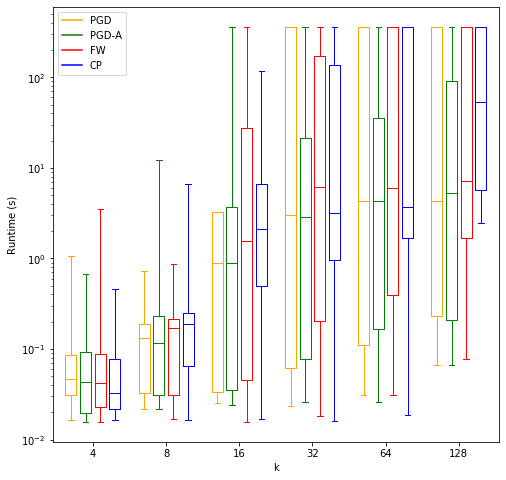}
        \caption{}\label{Boxplota}
    \end{subfigure}
    \begin{subfigure}{0.48\textwidth}
        \centering
        \includegraphics[scale = 0.3, trim = {0 40 0 65pt}, clip]{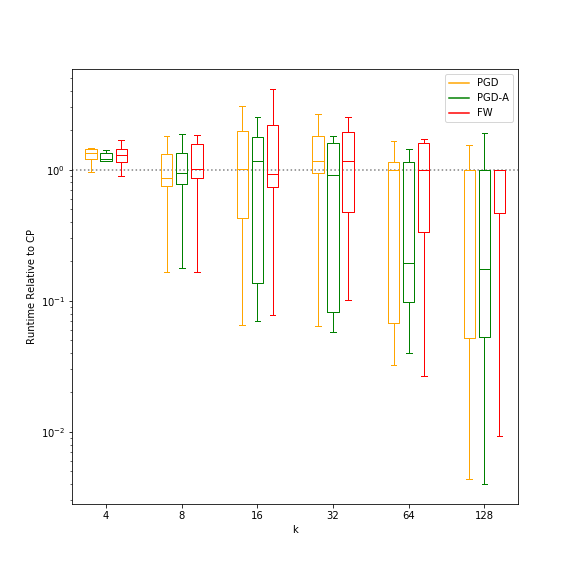}
        \caption{}\label{Boxplotb}
    \end{subfigure}
    \caption{Boxplots for each value of $k$ of: (a) the median runtimes across FOP instances, and (b) the ratio of median runtimes compared to CP median runtime for the 3 gap-gradient methods.}
    \label{fig:MIPLIBboxplots}
\end{figure}

Figure~\ref{fig:MIPLIBExp1} displays the median solution times of each novel algorithm in comparison to CP. PGD, PGD-A, and FW improved upon CP in 27, 30, and 26 of the 54 instances, respectively. For those instances where our methods did not improve upon CP, the performance was similar. Figure~\ref{fig:MIPLIBExp2} displays the median iteration count of each novel algorithm in comparison to CP. By design of the termination criterion, our gap-gradient methods take a minimum of 3 iterations before terminating, and for many instances where $k = 4$, the initialized value of $\balpha^{(0)} = \frac{1}{\mathcal{K}}\mathbf{1}$ is a loss function minimizer, and thus can be completed in one iteration by the CP method.

Figure~\ref{fig:MIPLIBboxplots} displays boxplots of the median solution time for each of the four methods and the improvement in median solution time over CP. The median solution time ($\pm$ standard deviation) as a function of $k$ for CP was 0.03s ($\pm$0.14), 0.19s ($\pm$2.01), 2.10s ($\pm$37.0), 3.17s ($\pm$137.2), 3.75s ($\pm$168.5), 53.1s ($\pm$160.0). The median solution time ($\pm$ standard deviation) for our best performing implementation (PGD-A) was 0.04s ($\pm$0.24), 0.12s ($\pm$3.82), 0.88s ($\pm$137.5), 2.91s ($\pm$111.1), 4.29s ($\pm$110.3), 5.29s ($\pm$128.1). The median improvement in median solution time as a function of $k$ was -0.01s (-33.3\%), 0.07s (36.8\%), 1.22s (58.1\%), 0.26s (8.2\%), -0.54s (-14.4\%), and 47.81s (90.0\%).

\section{Application to Political Gerrymandering} \label{s.res}

In this section, we present the FOP used for political districting (Section~\ref{FOP}), present custom solution approaches based on graph coarsening and ensemble methods (Section~\ref{Coarse}), and then evaluate our solution approaches using randomly generated problem instances (Section~\ref{Experiments}). In the following section (Section~\ref{sec:iowa}), we present a case study using real data from the State of Iowa.

\subsection{Forward optimization problem (FOP)}\label{FOP}

We use the political districting model from \cite{swamy2022multiobjective} as our FOP. The model is a mixed-integer linear program that determines the political districting for a given state. Let $G = \{V, E\}$ be a graph where each vertex $v \in V$ represents a census block or some larger tract of land within the state and edges $e \in E$ represent the adjacency of such areas. Each vertex $v$ is associated with a population $p_v \in V$, the number of democratic (republican) voters $p_v^D (p_v^R), v\in V$, and an area (in square feet) $a_v, v\in V$. We let $d_{i,j}, i,j\in V$ represent the Euclidean distance between vertex $i$ and $j$ (not necessarily adjacent). 

Let the decision variables $x_{ij}$ be a binary indicator of whether the vertex $i \in V$ is the center of a district that contains the vertex $j \in V$. The decision variable $f_{ijv}$ indicates a network flow between vertices $j$ and $v \in V$ where they both have district center $i$, which is used to maintain contiguity of modeled districts. The binary decision variable $z_{i}^D$ indicates if the district with center $i$ is won by the Democratic party, and $v_{ij}^D$ indicates whether or not $j$ is in a district with center $i$ which is won by the Democratic party. Finally, the decision variable $w_i$ indicates the number of wasted Democratic votes minus the number of wasted Republican votes in the election in the district with center $i$.

Intuitively, the model solution determines a partition of the graph $G$ into $L$ disjoint connected subgraphs, which denote the legislative districts. We provide the full FOP formulation in Appendix~\ref{app:FOP-formulation}. We define the objective of the FOP as a weighted sum of three components: a measure of population imbalance, a measure of compactness, and a measure of efficiency gap.

\subsubsection{Population imbalance ($\rho$).} We let $\rho$ represent the population imbalance measured as the largest relative deviation from the average district population exhibited by any district. 
$$\rho = \frac{\max_{i \in V} |\sum_{j \in V}p_i x_{ij} - \frac{\sum_{j \in V} p_j}{L}|}{\frac{\sum_{j \in V} p_j}{L}}$$

For example, in a state of population 1000, ten districts of population 100 would yield $\rho = 0$, while nine districts with population 101 and one of population 91 would yield $\rho = \frac{|91 - 100|}{100} = 0.09$. Unlike \cite{swamy2022multiobjective}, we include population imbalance in the objective rather than a constraint because the current legal text and judicial practice indicates that district populations should be \emph{as close to equal as possible}, clearly indicating it as an objective to be minimized.

\subsubsection{Compactness ($\sigma_A$).} We let $\sigma_A$ represent the compactness of a districting. Compactness is formulated as a measure of the p-median distance weighted by the area of each tract of land, divided by the area-weighted 1-median distance of the entire state (denoted by $M$). Using our notation, this can be written as:
$$\sigma_A = \frac{\sum_{i, j \in V}d_{ij} a_j x_{ij}}{M}.$$

\subsubsection{Efficiency gap ($\phi_{EG}$).} We let $\phi_{EG}$ denote the efficiency gap, measured as defined by \cite{stephanopoulos2015partisan} using a set of constraints formulated by \cite{swamy2022multiobjective}. The measurement of the efficiency gap is motivated by the construction and counting of \emph{wasted votes}, which contain both votes cast for a losing candidate and votes cast for a winning candidate beyond the necessary majority. The relative disparity in wasted votes distributed across two major parties as the result of district-based elections is measured as the efficiency gap, which can be calculated as:
$$\phi_{EG} = \frac{|\sum_{i \in V} w_i|}{\sum_{i \in V} p^D_i + p^R_i}$$

The construction of this metric directly targets districting tactics known as \emph{cracking} and \emph{packing}, which are alleged to be a major strategy of partisan gerrymandering. 

Since this problem is very computationally taxing to solve, \cite{swamy2022multiobjective} propose a graph coarsening approach to reduce the size of the problem formulation.

\subsection{Graph coarsening}\label{Coarse}

In this section, we present custom solution approaches for graph based inverse mixed integer optimization problems. Our approach uses a combination of graph coarsening with ideas from machine learning and data-driven inverse optimization. Although these methods were inspired by our application of political districting, they can be applied to any inverse optimization problem where the FOP can be represented as a graph-based problem.

\subsubsection{Maximal matching.}\label{altmaxmatch}

Basic graph coarsening is achieved by producing a random maximal matching of the edges of the graph and contracting all edges in the produced matching into vertices. A random maximal matching is produced by iterating through a random ordering of the edges, and selecting an edge to be contracted if both of its endpoints are not adjacent to any edges that have already been selected for contraction. All population and area data at the endpoints of a contracted edge is summed together to produce the data for the newly formed vertex, and the distances associated with the edges of newly formed vertices are calculated from landmass centroids derived from area-weighted averages of the two contracted endpoints. This process can be iteratively repeated as many times as desired to shrink the graph to a needed size. Overall, this process produces a smaller graph with the same underlying spatial structure but at a lower data resolution. This allows for a problem that is less computationally taxing, at a tradeoff of reduced accuracy of the final solution due to lower data resolution. 

\cite{swamy2022multiobjective} suggest that coarsening methods that choosing matchings that prioritize merging vertices with a low combined population preserve solutions with lower population imbalances, and lower objectives values for many subobjectives, with the same improvement on computational efficiency as random maximal matchings. As such, their method is deterministic, producing an ordering of edges by directly ranking the combined populations.

We propose a variation that still involves random choices of edges (which we exploit in the next section), while gaining the advantage of coarsenings that preserve better solutions. As such, our method creates an ordering of edges by sampling an exponential random variable for each edge with a mean of the combined population divided by twice the mean vertex population, and ordering by the edges by the values of the random samples. We note that if all edges have identical combined populations, then this method produces the same distribution of outputs as a random maximal matching method. 

\subsubsection{Ensemble-based coarsening.}
We leverage the idea of \emph{ensembles} from machine learning to produce and exploit multiple coarsenings of the same graph \citep{breiman1996bagging}. Since both the basic maximial matching and our proposed variation are non-deterministic, multiple (different) coarsenings of the same graph can be produced with the same method. While each coarsening loses some resolution, elements lost in one coarsening may be preserved in another. An ensemble-based approach may allow for the risks of coarsening to be reduced in the aggregate, while potentially preserving the computational savings of coarsening. Let FOP$_1$, FOP$_2$ ... FOP$_n$ denote the forward optimization problems generated for an ensemble of $n$ different coarsenings of a political districting FOP, and let $\mathcal{S}_1, \mathcal{S}_2 \dots \mathcal{S}_n$ denote their respective feasible regions. We present two options for how an ensemble of coarsenings can be used to produce heuristic solutions to the inverse problem:

\begin{enumerate}
    \item \textbf{Solve each instance independently.} In this approach, we solve an inverse model independently on each coarsened graph, generating the inverse solution $\balpha^*_i$ from each individual inverse formulation GIO$_\text{ABS}(\hat\by, \mathcal{S}_i)$. The average of the outputs, $\frac{1}{n} \underset{i \in 1 \dots n}{\sum} \balpha^*_i$, can then be used as an approximation of the true $\balpha^*$. Alternatively, the set of inverse solutions can be analyzed as a distribution. For example, the convex hull of the set of outputs can be interpreted as a polytope of potential inverse solutions.
    \item \textbf{Solve a multi-point formulation.} In this approach, we solve a single multi-point inverse optimization formulation produced by creating a master problem with constraints that are derived from multiple coarsened graphs at once. A single slack variable describes the optimality gap across all constraints, thus minimizing the maximum optimality gap across all coarsenings.
\end{enumerate}

For the remainder of this section, we focus on the second approach, which can be formalized as the following optimization problem:
\begin{align*}
    \minimize_{\bm{\alpha},\; \xi_{\text{ENS}}} \quad & \xi_{\text{ENS}}  \tag{MultiGIO$_{\text{ABS}}$MinMax}\\
    \text{subject to } \quad & \bm{\alpha}^\intercal \bC \bm{y^0} \leq \bm{\alpha}^\intercal \bC \bm{y} + \xi_{\text{ENS}}, \quad \forall \: \bm{y} \in \mathcal{S}_1, \mathcal{S}_2, \dots \mathcal{S}_n, \\
    &||\bm{\alpha}||_1 = 1, \\
    &\bm{\alpha} \geq \bm{0}, \:\: \xi_{\text{ENS}} \geq 0. 
\end{align*}

For MultiGIO$_{\text{ABS}}$MinMax, we can show that as as the size of the ensemble approaches infinity, the likelihood of our method finding the same solution as the full resolution graph approaches $1$.

\begin{theorem}\label{thm:ensemble_converge}
    Given an ensemble of $n$ independently sampled coarsenings of $G$ down to $v \geq L$ vertices, such that any coarsened graph with $v$ vertices has a non-zero probability of being sampled, as $n \rightarrow \infty$, the probability of MultiGIO$_{\text{ABS}}$MinMax yielding $\balpha^*=$ approaches 1.
\end{theorem}

Proof of Theorem~\ref{thm:ensemble_converge} is provided in Appendix~\ref{ECproofensembleconverge}. In the context of implementing gap-gradient methods for the multipoint inverse formulation, the coarsening that returns the highest optimality gap at the hypothesis $\balpha^{(k)}$ determines the next step in the descent process.

\begin{proposition}\label{prop:multipoint_gradient}
    Let FOP$_1$, FOP$_2$ ... FOP$_n$ denote the forward optimization problems generated for an ensemble of $n$ different coarsenings of a political districting FOP, and let $\xi_\text{ENS}(\balpha)$ denote the gap function corresponding to the multipoint inverse formulation MultiGIO$_\text{ABS}$MinMax generated by this ensemble of coarsenings. Let $\balpha^{(k)}$ be the hypothesis cost vector at iteration $k$ of a solution method for MultiGIO$_\text{ABS}$MinMax, and let $y^{(k)}_i$ denote the optimal solution of FOP$_i(\balpha^{(k)})$. At any $\balpha^{(k)}$, the subgradient of $\xi_\text{ENS}(\balpha)$ is equal to $\bC (\hat\by - \underset{\by \in \by^{(k)}_1 \dots \by^{(k)}_n}{\argmin} (\balpha^{(k)\intercal}\bC\by))$.
\end{proposition}

Proof of Proposition~\ref{prop:multipoint_gradient} is provided in Appendix~\ref{ECproofmultipointgradient}.

Since the loss function $\xi_\text{ENS}(\balpha)$ satisfies all the properties listed in Remark~\ref{rem:graddesc}, the gap-gradient methods described in Section~\ref{sec:gap-gradient} may be used to find solutions to MultiGIO$_\text{ABS}$MinMax.

\subsubsection{Stochastic descent with coarsened ensembles.}

In this section, we demonstrate how stochastic gradient descent can be used with coarsening to solve generalized inverse optimization problems with graph-based FOPs. Instead of generating an ensemble of coarsenings to create a lower approximation $\xi_\text{ENS}$ of $\xi_\text{ABS}$, and then deterministically minimizing the approximation, one can instead use ensembles of coarsenings to produce estimators of subgradients of $\xi_\text{ABS}$ in the application of a stochastic gap-gradient method to an inverse model formulated for the full sized graph. Similar to the multipoint formulation method, our gradient approximation is calculated as $\nabla\xi_\text{ABS}(\balpha^{(k)}) \approx \bC (\hat\by - \underset{\by \in \by^{(k)}_1 \dots \by^{(k)}_n}{\argmin} (\balpha^{(k)\intercal}\bC\by))$. However, in this application, an ensemble of coarsenings FOP$_1^{(k)}$,  FOP$_2^{(k)}$\dots FOP$_n^{(k)}$ is independently randomly sampled during each iteration $k$ of the subgradient method. Here the use of coarsenings and ensembles is analogous to the process of mini-batching in stochastic gradient descent for training neural networks, where the true loss function gradient is estimated by calculating the loss function with respect to a randomly sampled subset of the the full set of available samples \citep{schmidt2019cpsc}. In neural networks, the samples come from the set of training observations; in the inverse graph partitioning case, `samples' come from the set of possible partitions $\mathcal{S}$, and each coarsened graph contains contains a subset of such samples.

It is not obvious that solving the FOP on coarsened graphs or ensembles thereof yields an unbiased estimator of $\nabla\xi_\text{ABS}$. However, Theorem~\ref{thm:ensemble_converge} does show that approximating the gradient of $\xi_\text{ABS}$ with $\bC (\hat\by - \underset{\by \in \by^{(k)}_1 \dots \by^{(k)}_n}{\argmin} (\balpha^{(k)\intercal}\bC\by))$ constitutes a \textit{consistent and asymptotically unbiased estimator} of $\nabla \xi_\text{ABS}$, as defined by \cite{chen2019stochastic}. Under the conditions of a convex but not strongly convex gap function (such as absolute sub-optimality), when using a subgradient estimation method where expected estimation error is inversely proportional to $\sqrt{k}$, stochastic projected gradient descent observes a convergence rate of $O(\frac{1}{\sqrt{k}})$ towards the gap function minimum value. To obtain this bound for the coarsening application, the subgradient estimation must involve progressively larger ensembles with each iteration. Without progressively increasing the ensemble size, we can expect our method to converge to the neighborhood of the minimizer, but is not guaranteed to find the exact minimizer.

Algorithm~\ref{stochasticggproj} presents the structure of an implementation of stochastic subgradient estimation within a projected gradient descent approach to loss function minimization, where $n_k$ indicates a rule for selecting a coarsening ensemble size at a given iteration $k$, and $K$ denotes a maximum number of iterations. At the end of $K$ iterations, a weighted averaging of the values of $\balpha^{(k)}$ is returned.
We note that the termination criterion utilized in Section~\ref{sec:gap-gradient} are much less likely to be triggered at any given iteration. 
As such, in our evaluation of stochastic subgradient methods, we examine how closely the descent methods approach the true loss function minimum and minimizer, rather then comparing the time until termination yielding an exact minimizer, as in the previous computational evaluations.

\subsubsection{Boosted ensembles of coarsenings.} \label{sec:ensemble_boosting}
To accentuate the diversity of the ensemble of coarsened graphs, we implement a method of sequential coarsening inspired by boosting ensemble methods from machine learning \citep{schapire1999brief}. In an unboosted coarsening implementation, the random maximal matching of a graph is determined by creating a uniformly randomly ordered list of the edges, proceeding through the list in order, and contracting each edge if and only if both endpoints are not adjacent to an edge that has already been contracted. The random shuffling that determines each coarsening is independent for each element of the ensemble. In a \emph{boosted implementation}, the randomly ordered shuffling of the edges is done in a weighted fashion such that an edge with a higher weight is more likely to be later in the ordering. In coarsening the first element of the ensemble, either the initial edge weights are uniform, or weighted by population as detailed in Section~\ref{altmaxmatch}. If an edge $e \in E$ is contracted when creating the $i^{\text{th}}$ element of the ensemble, then for the $(i+1)^{\text{th}}$ element, the weight of edge $e$ is multiplied by a factor $\eta > 1$, and its weight is multiplied by $\frac{1}{\eta}$ otherwise. This ensures that every ensemble member is less likely to have certain edges contracted that were frequently contracted in previous coarsenings, discouraging highly correlated coarsenings in the ensemble. This process is detailed in EC\ref{sec:EC_algstruc} (Algorithm \ref{alg:boosting}). In the case of ensembles of graphs coarsened multiple times over, this process can be repeated at each tier of coarsening in a tiered tree-like fashion until the desired depth of coarsening is reached.


\subsection{Computational Experiments}\label{Experiments}

In this section, we evaluate our solution approaches
using randomly generated problem instances representative of political districting.

\subsubsection{Experimental setup.} \label{districtExperimentSetup} 
To evaluate our models, we use 8 simulated graphs/states of size $|G| = 20$. See \ref{sec:exp} for details on how we simulate a state. For each graph/state, we create 5 uniformly randomly sampled multi-objectives weightings, and solve the corresponding FOP to create 5 inverse inputs for each simulated state. In total, we have 40 simulated IO instances. For each instance, we perturb the value of each subobjective by adding a random value in the range [0.0375, 0.0625]; this allows us to find an interior point that is likely to have a unique gap function minimizer. We solve each instance to optimality without any coarsening as a benchmark. 

We conduct three experiments. 
First, we evaluate the impact of coarsening on the accuracy of inverse optimization solutions. We evaluate two coarsening methods: contracting edges selected from a random maximal matching of vertices and contracting edges selected from a population-weighted maximal matching. For each method, we conduct both one and two rounds of coarsening for each of the 40 problem instances. We then solve an IO problem for each instance using the PGD-A method to obtain a heuristic solution $\balpha^h$. 

Second, we evaluate the efficacy of the multi-point ensemble formulation as a heuristic for minimizing the gap function. We generate coarsened graphs with one round of random maximal matching, and we create ensembles of size 1, 4, 16, and 64, using both the independently random coarsening and the boosted coarsening techniques described in Section~\ref{sec:ensemble_boosting} with learning rate $\eta = 1.5$. For each ensemble, we apply the multi-point inverse formulation MultiGIO$_{\text{abs}}$MinMax on the ensemble to achieve the heuristic solution $\balpha^h$. 

Third, we evaluate the performance of stochastic descent methods on estimating the gap function minimizer. We apply stochastic descent methods with unboosted ensembles of sizes 16 and 64, and an increasing ensemble rule $n=k$. For the descent methods, a step size is automatically chosen such that the first step has length 0.1, with following step sizes decreasing at a rate of $\frac{1}{\sqrt{k}}$, and a momentum coefficient of 0.1. We run each descent algorithm for a total of 12 iterations and return a heuristic solution $\balpha^h$ at iterations 0, 4, 8, and 12. 

For each experiment, we record the total solution time, the gap function value of the returned output cost vector minus the true gap function minimum value (i.e., $\xi_\text{ABS}(\balpha^h) - \xi_\text{ABS}(\balpha^*)$), and the Euclidean distance between the returned cost vector and the optimal cost vector (i.e., $||\balpha^h - \balpha^*||_2$).

\subsubsection{Experimental results.}

\begin{figure}
    \centering
    \begin{subfigure}{0.45\textwidth}
        \centering\includegraphics[scale=0.4,trim={0 0 0 0pt},clip]{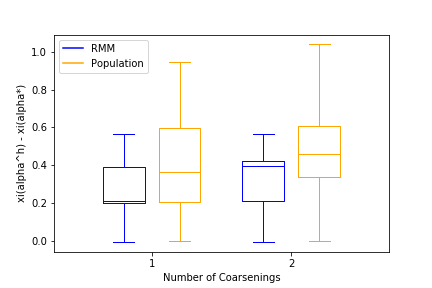}
        \caption{}\label{Exp3Ba}
    \end{subfigure}
    \begin{subfigure}{0.45\textwidth}
        \centering\includegraphics[scale=0.4,trim={0 0 0 20pt},clip]{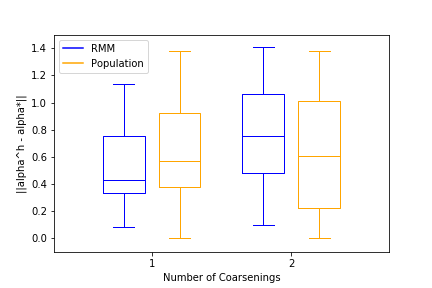}
        \caption{}\label{Exp3Bb}
    \end{subfigure}
    \caption{A comparison of (a) the difference between retrieved and true optimality gap, and (b) the distance from the original objective for different numbers of coarsenings.}
    \label{fig:Exp3B}
\end{figure}

Figure~\ref{Exp3Ba} displays the gap function value of the returned output cost vector minus the true gap function minimum value (i.e., $\xi_\text{ABS}(\balpha^h) - \xi_\text{ABS}(\balpha^*)$). For coarsening by random maximal matching, the median (standard deviation) difference was 0.21 (0.14) and 0.39 (0.16) for 1 and 2 rounds of coarsening, respectively. For coarsening by population-weighted matching, the median (standard deviation) difference was 0.37 (0.26) and 0.46 (0.27) for 1 and 2 rounds of coarsening, respectively. Figure~\ref{Exp3Bb} displays the Euclidean distance between the returned cost vector and the optimal cost vector (i.e., $||\balpha^h - \balpha^*||_2$). For coarsening by random maximal matching, median (standard deviation) distance was 0.43 (0.31) and 0.75 (0.33) for 1 and 2 rounds of coarsening, respectively. For coarsening by population-weighted matching, median (standard deviation) distance was 0.57 (0.33) and 0.60 (0.43) for 1 and 2 rounds of coarsening, respectively. Overall, random maximal matching produced solutions that were closer to the optimal solution as compared to population-weighted matching.

\begin{figure}
    \begin{subfigure}{0.3\textwidth}
        \includegraphics[scale=0.4,trim={0 0 0 0},clip]{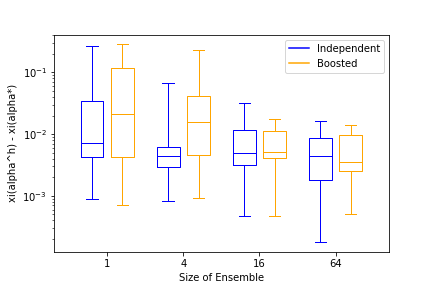}
        \caption{}\label{Exp4Aa}
    \end{subfigure}\hspace*{2em}
    \begin{subfigure}{0.3\textwidth}
        \includegraphics[scale=0.4,trim={0 0 0 0},clip]{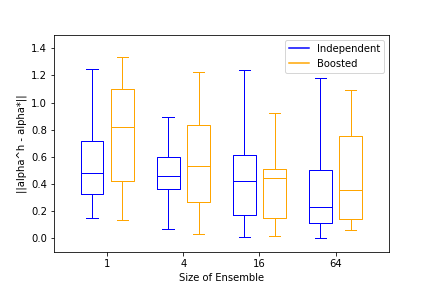}
        \caption{}\label{Exp4Ab}
    \end{subfigure}\hspace*{2em}
    \begin{subfigure}{0.3\textwidth}
        \includegraphics[scale=0.4]{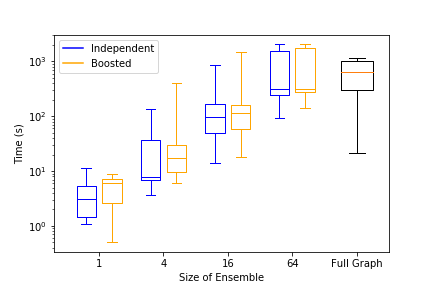}
        \caption{}\label{Exp4Ac}
    \end{subfigure}
    \caption{A comparison of (a) the difference between retrieved and true optimality gap, (b) the distance from the original objective, an (c) the solution time for various ensemble sizes.}
    \label{fig:Exp4A}
\end{figure}

Figure~\ref{Exp4Aa} displays $\xi_\text{ABS}(\balpha^h) - \xi_\text{ABS}(\balpha^*)$. For unboosted ensembles, the median (standard deviation) difference was 0.0071 (0.075), 0.0045 (0.014), 0.0050 (0.0079), and 0.0045 (0.0046), for ensembles of size 1, 4, 16, and 64, respectively. For boosted ensembles, the median (standard deviation) difference was 0.021 (0.092), 0.016 (0.048), 0.0051 (0.0053), and 0.0036 (0.0043), for ensembles of size 1, 4, 16, and 64, respectively. Overall, an unboosted ensemble of size 64 was able to reduce the median and standard deviation of $\xi_\text{ABS}(\balpha^h) - \xi_\text{ABS}(\balpha^*)$ by 0.0026 (37\%) and 0.704 (93\%) over an ensemble of size 1, respectively. Figure~\ref{Exp4Ab} displays $||\balpha^h - \balpha^*||_2$. For unboosted ensembles, the median (standard deviation) distance was 0.48 (0.28), 0.46 (0.20), 0.42 (0.31), and 0.23( 0.30), for ensembles of size 1, 4, 16, and 64, respectively. For boosted ensembles, the median (standard deviation) distance was 0.82 (0.36), 0.53 (0.36), 0.45 (0.26), and 0.35 (0.33), for ensembles of size 1, 4, 16, and 64, respectively. Overall, an unboosted ensemble of size 64 was able to reduce the median distance by 0.25 (52\%) over an ensemble of size 1. Figure~\ref{Exp4Ac} displays the solution time; an unboosted ensemble of size 64 was able to reduce the median and standard deviation of the solution time of the full graph by 341s(52\%) and -396s (-10\%), respectively.

\begin{figure}
    \begin{subfigure}{0.3\textwidth}
        \includegraphics[scale=0.4,trim={0 0 0 0},clip]{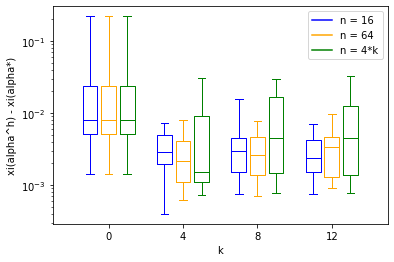}
        \caption{}\label{stochexp1}
    \end{subfigure}\hspace*{2em}
    \begin{subfigure}{0.3\textwidth}
        \includegraphics[scale=0.4,trim={0 0 0 0},clip]{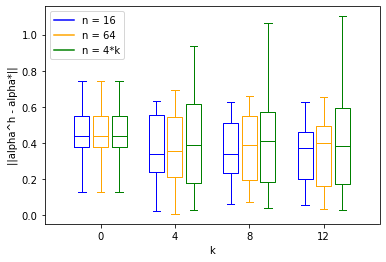}
        \caption{}\label{stochexp2}
    \end{subfigure}\hspace*{2em}
    \begin{subfigure}{0.3\textwidth}
        \includegraphics[scale=0.4]{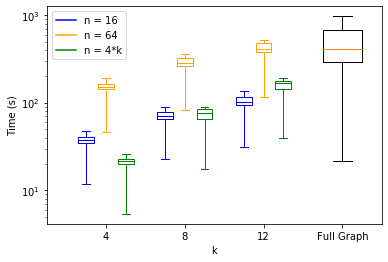}
        \caption{}\label{stochexp3}
    \end{subfigure}
    \caption{A comparison of (a) the difference between retrieved and true optimality gap, (b) the distance from the original objective, (c) the running time.}
    \label{fig:StochExp}
\end{figure}

Figure~\ref{stochexp1} displays $\xi_\text{ABS}(\balpha^h) - \xi_\text{ABS}(\balpha^*)$. For stochastic descent, the median (standard deviation) difference after 12 iterations was 0.0031 (0.0047), 0.0037 (0.0036), and 0.0056 (0.0083), for ensembles of size 16, 64, and increasing size, respectively. Figure~\ref{stochexp2} displays $||\balpha^h - \balpha^*||_2$. The median (standard deviation) difference after 12 iterations was 0.37 (0.18), 0.40 (0.19), and 0.38 (0.31), for ensembles of size 16, 64, and increasing size, respectively. Figure~\ref{stochexp3} displays the solution time; 12 iterations of stochastic descent with an ensemble size of 64 was able to reduce the median and standard deviation solution time of the full graph by 55s (11\%) and 211s (71\%), respectively.

Overall, these results suggest that a multi-point formulation with an ensemble of 64 unboosted graphs coarsened by random maximal matching performs best in providing accurate approximations, while still providing decreased solution times.

\section{Case Study: the State of Iowa}\label{sec:iowa}

To demonstrate the application of inverse optimization to political gerrymandering, we applied our methods to the current congressional districts for the state of Iowa. 

\subsection{Context}

Iowa's 2022 congressional districts were designed by an independent districting commission. The first set of districts submitted by the commission were rejected by the state legislature. A second set of districts returned by the commission was subsequently approved by the state legislature. 

Congressional districting in Iowa is unique in that the state constitution requires that counties are not split by district borders. The Iowa constitution also provides specific guidelines on how district compactness is to be measured (Iowa Code section 42.2). Two possible metrics are defined; \emph{length-width compactness}, which is calculated as the absolute value of the difference in a district's east-west distance and its north-south distance, and \emph{perimeter compactness}, which is calculated as the perimeter of a district. We modify our FOP to measure compactness as the sum of the perimeter distance of each district, divided by the state perimeter times the number of districts (to scale the metric to the same range as the other sub-objectives). Perimeter compactness was chosen over length-width compactness as there exist numerous examples of possible districts with near-0 length-width compactness that by most common-sense approaches may be considered very non-compact, and while perimeter compactness may be susceptible to `coastline paradoxes' that can make perimeter measurements not correspond to the broader shape of a district, the majority of Iowa's counties are rectilinear in their boundaries, making this issue mostly irrelevant. 

Let $\sigma_P$ denote the perimeter compactness, let $q_{ij}$ be binary variable indicating if adjacent vertices $i$ and $j$ are in different districts, let $b_{ij}$ represent the length of the border shared by counties $i$ and $j$, and let $M_p$ represent the perimeter distance of the state of Iowa. We represent the metric using the following constraints:
\begin{subequations} \label{FP}
\begin{align}
    q_{ij} &\geq x_{ki} - x_{kj}, &\forall i, j, k &\in V, b_{ij} \neq 0, \tag{2t}\\
    \sigma_P &= \frac{(\sum_{i, j \in V}q_{ij} b_{ij}) + M_p}{L M_p}, \tag{2u}\\
     q_{ij} &\in \{0,1\},  &\forall i, j &\in V, b_{ij} \neq 0. \tag{2v}
\end{align}
\end{subequations}

\subsection{Data}

\begin{figure}
    \centering
    \includegraphics[scale=0.3]{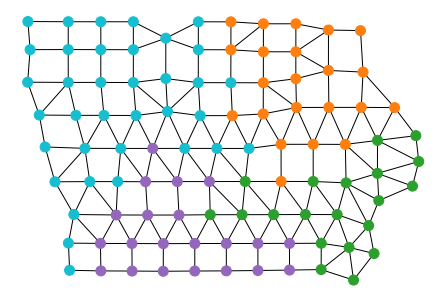}
    \includegraphics[scale=0.3]{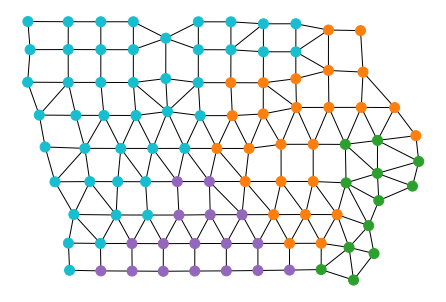}
    \includegraphics[scale=0.3]{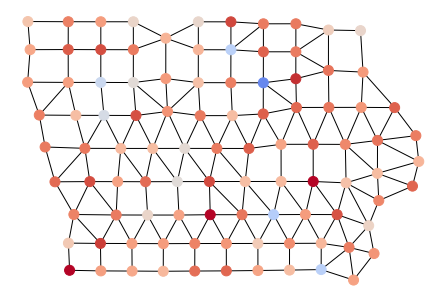}
    \caption{The full-sized graph of Iowa Counties, with color and shape denoting (a) the enacted 2022 districts, (b) the initially proposed and rejected 2022 districts, and (c) by relative partisan slant of the population in 2020 statewide elections}
    \label{fig:stategraph}
\end{figure}

State data for Iowa districts and electoral data were obtained from the Metric Geometry and Gerrymandering Group and the ALARM Project \citep{mccartan2022simulated}. The full-sized graph of Iowa comprises 99 counties, partitioned into 4 districts.  Figure \ref{fig:stategraph} shows the full-sized graph of the state. The first subfigure displays the current assignment of counties into 4 districts (4 colours), the second subfigure displays the initially proposed 4 distrcts, and the final subfigure displays the relative partisan slant of voters in 2020 statewide elections, with blue indicating more Democrats and red indicating more Republicans.

\begin{figure}
    \centering
    \includegraphics[scale=0.4,trim={0 0 0 20pt},clip]{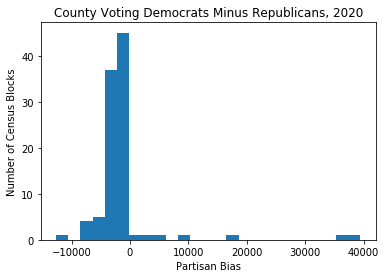}
    \includegraphics[scale=0.4,trim={0 0 0 20pt},clip]{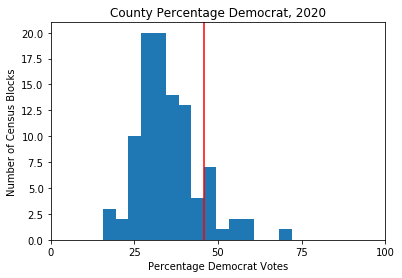}
    \caption{Distribution of (a) number of Democrat votes minus Republican votes and (b) percentage of votes towards Democrats by county}
    \label{fig:EDA1}
\end{figure}

Figure \ref{fig:EDA1} (a) and (b) display the distribution of the difference between Democrat and Republican votes across the Iowa counties in the 2020 presidential election (the positive side of the x-axis represents more Democrat votes) by absolute difference and per capita percentages, respectively. The vertical line in subfigure (b) indicates the partisan slant of the entire whole.  Democratic votes make up 46\% of all votes cast, but 93 out of the 99 of counties have republican majorities, with the democratic majorities being by larger margins in fewer counties. As a result the counties themselves display a noticeable efficiency gap, but rather than attributing this to gerrymandering of county borders, it may be possible that socio-geographical phenomena (e.g., urban/rural divides) make it such that Democratic and Republican populations are generally distributed as such in the Iowa political geography, regardless of borders. This illustrates the importance of understanding partisan fairness metrics such as efficiency gap in terms of what is actually feasible in the given political geography, and in relation to the existence of other competing metrics.

The sub-objective values of our inverse inputs are detailed in Table \ref{tab:IowaInvInput}. Population data is drawn from the 2020 Census and political data from the 2020 presidential election. We note that in both districtings, the evaluation of the efficiency gap is in favor of Republican over-representation. We also supply district-by-district measurements of the metrics that are used to calculate inverse inputs in Tables \ref{tab:IowaRejectedDistrictMetrics} and \ref{tab:IowaAcceptedDistrictMetrics} for the rejected and accepted districting plans, respectively. For context we note that the ideal district population is 797,592 and the Iowa state perimeter is 1151 miles.

\begin{table}[]
    \centering
    \begin{tabular}{c|c|c}
         \textbf{Metric} & \textbf{Rejected Value} & \textbf{Accepted Value}\\
         \hline 
         Perimeter Compactness $(\sigma_P)$ & 0.5773 & 0.6116 \\
         Population Imbalance $(\rho)$ & $7.8674 * 10^{-5}$ & $6.6137 * 10^{-5}$ \\
         Efficiency Gap $(\phi_{EG})$& 0.0882 & 0.4163\\
    \end{tabular}
    \caption{Iowa Case Study Inverse Inputs}
    \label{tab:IowaInvInput}
\end{table}

\begin{table}[]
    \centering
    \begin{tabular}{c|c|c|c|c}
         \textbf{District} & \textbf{Population} & \textbf{Population Deviation} & \textbf{Perimeter} & \textbf{\% Democrat}\\
         \hline 
         1 & 797,655 & 0.008\% & 503.03 & 54.6\% \\
         2 & 797,556 & 0.005\% & 786.27 & 44.9\% \\
         3 & 797,584 & 0.001\% & 515.54 & 49.9\% \\
         4 & 797,574 & 0.002\% & 967.18 & 34.9\%\\
    \end{tabular}
    \caption{Iowa Rejected Individual District Metrics}
    \label{tab:IowaRejectedDistrictMetrics}
\end{table}

\begin{table}[]
    \centering
    \begin{tabular}{c|c|c|c|c}
         \textbf{District} & \textbf{Population} & \textbf{Population Deviation} & \textbf{Perimeter} & \textbf{\% Democrat}\\
         \hline 
         1 & 797,584 & 0.001\% & 696.48 & 49.0\% \\
         2 & 797,589 & 0.00037\% & 624.17 & 48.1\% \\
         3 & 797,551 & 0.0051\% & 619.72 & 49.8\% \\
         4 & 797,645 & 0.0066\% & 997.50 & 37.4\%\\
    \end{tabular}
    \caption{Iowa Accepted Individual District Metrics}
    \label{tab:IowaAcceptedDistrictMetrics}
\end{table}

\subsection{Inverse Optimization Analysis}

Although Iowa is computationally easier than other larger states, coarsening and ensemble methods are still needed for tractability. We use our heavy ball descent method, in conjunction with an ensemble of 64 graphs, each coarsened from the original state data three times. For our preliminary inverse analysis of the state, we produce an ensemble of 64 coarsenings of the state using 3 rounds of random maximal matching. Over this ensemble of coarsenings, we solve the MultiGIO$_{\text{abs}}$MinMax multipoint formulation, using PGD-A  for finding the optimality gap minimizer with the same parameters described in Section~\ref{districtExperimentSetup}.

\begin{table}[]
    \centering
    \begin{tabular}{c|c|c}
         \textbf{Metric} & \textbf{Analysis 1 Objectives} & \textbf{Analysis 2 Objectives}\\
         \hline 
         Perimeter Compactness & 0.055 & 0 \\
         Population Imbalance &  0.933 & 1 \\
         Efficiency Gap & 0.011 & 0 \\
         \hline Reported Optimality Gap & -0.0014 & -0.0014 \\
    \end{tabular}
    \caption{Iowa 2022 Rejected Districts Inverse Output}
    \label{tab:IowaRejInvOutput}
\end{table}
\begin{table}[]
    \centering
    \begin{tabular}{c|c|c}
         \textbf{Metric} & \textbf{Analysis 1} & \textbf{Analysis 2}\\
         \hline 
         Perimeter Compactness & 0 & 0 \\
         Population Imbalance & 1 & 1 \\
         Efficiency Gap & 0  & 0 \\
         \hline Reported Optimality Gap & 0.0772 & 0.0042\\
    \end{tabular}
    \caption{Iowa 2022 Accepted Districts Inverse Output}
    \label{tab:IowaAccInvOutput}
\end{table}

 We conduct two inverse optimization analyses for the accepted and rejected districtings. In Analysis 1, we minimize the optimality gap for a objective weighting that is a combination of minimizing perimeter compactness ($\sigma_P$), population imbalance ($\rho$), and efficiency gap ($\phi_{EG}$). In Analysis 2, we  reformuate the efficiency gap metric so that the objective component represents \textit{maximizing} efficiency gap in favor of the Republican party (which both districtings favor). Since the original efficiency gap objective is formulated as an absolute value of a quantity $\frac{\sum_{i \in V} w_i}{\sum_{i \in V}(p^D_i + p^R_i)}$ that is positive when it benefits Democrats and negative when it benefits Republicans (see equations (2p) and (2q)), we substitute $\phi_{EG}$ in the objective with this `directional' quantity, and for our inverse input, the value of this quantity is equal to the negative efficiency gap. Minimizing this quantity is equivalent to maximizing the amount that the efficiency gap favors Republicans. We note that the linear dependence of the minimizing and maximizing $\phi_{EG}$ objectives makes it impossible for both subobjectives to be included in the same inverse analysis as the FOP objective feasible space might not have a fully dimensional feasible region. Thus by doing side-by-side analyses with both objectives and comparing the found optimality gaps, we can assess these objectives in comparison.
 
The results from both analyses are detailed in Tables \ref{tab:IowaRejInvOutput} and \ref{tab:IowaAccInvOutput} for the rejected and accepted districting plans, respectively. All analyses indicate that both districting plans place nearly all emphasis on minimizing population imbalance, with only the rejected districting placing some objective weight on minimizing perimeter compactness and partisan efficiency gap. 

The interpretation of our inverse analyses of these two districtings is that the accepted districting prioritizes minimizing population imbalance with. This supports the hypothesis that the currently enacted Iowa districts are not intentionally gerrymandered to benefit the Republican party. However, the extra priority placed on minimizing the population imbalance between the rejected and accepted plans in order to reduce imbalance by 16\% also increases to the perimeter and efficiency gap metrics by 6\% and 372\% respectively. Thus, this may still raise questions and discussion as to whether sole focus on reducing population imbalance is necessary when it entails such costs to other democratic values, even if there is not explicit intent to gerrymander districts.

We note that alternative quantitative analysis methods in the literature often model equal population considerations in districtings through constraints that put an upper limit on the population imbalance in modeled districts; the ALARM Project \citep{mccartan2022simulated} simulates possible districts with an upper limit of 0.5\% population imbalance, and \cite{swamy2022multiobjective} utilize an upper bound of 4.8\% in analyzing possible district plans for Wisconsin. We note that our application of generalized inverse optimization methods allows for the addition of population imbalance as a metric to be considered as optimized in itself (whereas \cite{swamy2022multiobjective}'s method of plotting Pareto frontiers is not practical for considering more than two objectives at once). As such, our methods are able to bring new considerations that indicate that with this objective considered, Iowa's districts may not be the result of prioritizing Republican advantage, whereas the ALARM Project's simulation-based analysis that permits larger population imbalances indicates that the efficiency gap of the employed districting is anomalous to the point of being indicative of gerrymandering.

We note that the analysis of the rejected districting reports optimality gaps that are negative, which indicates that our inverse input objective values cannot be Pareto-dominated by any solution that is feasible for some element in our ensemble of coarsened graphs. Since the Iowa districting input is in fact a feasible districting, it should be possible that a more comprehensive inverse modeling analysis, either utilizing the full-sized graph of Iowa without any heuristic methods, or an analysis with larger ensembles and/or graphs that are coarsened to a lesser extent, could retrieve a solution with a non-negative reported optimality gap. However, if the Pareto frontier of the feasible region covered by our ensemble is reasonably similar
to the geometry of the true forward problem Pareto frontier, then the results are still interpretable.

\section{Conclusion} \label{s.conclusion}

In summary, this paper makes three contributions. First, we propose a new approach for solving generalized mixed-integer IO problems based on sub-gradient methods. Second, we develop custom heuristic methods for graph-based inverse problems using a combination of graph coarsening and ensemble methods. Third, we propose a new application domain – quantitatively identifying gerrymandering – for generalized inverse integer optimization. We argue that IO can produce more nuanced data-driven arguments that proposed districtings should be considered gerrymandered.

\ACKNOWLEDGMENT{The authors thank Professor Ian Zhu for providing insightful feedback and suggestions to improve this work.}

\bibliographystyle{ormsv080}
\bibliography{references}




\newpage
\ECSwitch
\ECHead{Electronic Companion}

\section{Proofs} \label{ECproofs}

\subsection{Proof of proposition~\ref{lemma:MPsolve}} \label{ECproofMPsolve}

\proof{Proof of Proposition~\ref{lemma:MPsolve}:} Let $\xi^k_{ABS}(\balpha)$ denote the gap function corresponding to the inverse optimization model GIO$_{\text{ABS}}(\hat\by, \mathcal{S}^k)$. We note that that GIO$_{\text{ABS}}(\hat\by, \mathcal{S}^k)$ is a relaxation of GIO$_{\text{ABS}}(\hat\by, \mathcal{S})$, and as such $\xi^k_{ABS}(\balpha) \leq \xi_{ABS}(\balpha)$. We prove the claim by showing that if $\bC \hat{\by} - \xi_{\text{ABS}}(\balpha^*)\mathbf{1} \in \text{conv}(B^k)$ then; (i) $\xi^k_{ABS}(\balpha^*) = \xi_{ABS}(\balpha)$, and (ii) $\xi^k_{ABS}(\balpha^*) < \xi^k_{ABS}(\balpha), \: \forall \balpha \in \mathcal{A} \setminus \balpha^*$.

\textit{Proof of (i)}: Since $\bC \hat{\by} - \xi_{\text{ABS}}(\balpha^*)\mathbf{1} \in \text{conv}(B^k)$, it is true that $\bC \hat{\by} - \xi_{\text{ABS}}(\balpha^*)\mathbf{1} = \lambda_0 \bC \by^{(0)} + \lambda_1 \bC \by^{(1)} + \dots + \lambda_k \bC \by^{(k)}$ for some $\mathbf{\lambda} \in [0,1]^k$, $||\mathbf{\lambda}|| = 1$. Thus we have,
\begin{align*}
    \xi_{\text{ABS}}(\balpha^*) &= \balpha^{*\intercal} \bC \hat{\by} - \underset{\by \in \mathcal{S}}{\min} \: \balpha^{*\intercal} \bC \by \\
    \balpha^{*\intercal} \bC \hat{\by} - \balpha^{*\intercal}(\bC \hat{\by} - \xi_{\text{ABS}}(\balpha^*)\mathbf{1}) &= \balpha^{*\intercal} \bC \hat{\by} - \underset{\by \in \mathcal{S}}{\min} \: \balpha^{*\intercal} \bC \by \\
    \balpha^{*\intercal} \bC \hat{\by} - \balpha^{*\intercal}(\lambda_0 \bC \by^{(0)} + \lambda_1 \bC \by^{(1)} + \dots + \lambda_k \bC \by^{(k)}) &= \balpha^{*\intercal} \bC \hat{\by} - \underset{\by \in \mathcal{S}}{\min} \: \balpha^{*\intercal} \bC \by \\
    \balpha^{*\intercal}(\lambda_0 \bC \by^{(0)} + \lambda_1 \bC \by^{(1)} + \dots + \lambda_k \bC \by^{(k)}) &= \underset{\by \in \mathcal{S}}{\min} \: \balpha^{*\intercal} \bC \by
\end{align*}
This implies that each $\by^{(i)} \in S^k$ for which $\lambda_i > 0$ is a minimizer of $\underset{\by \in \mathcal{S}}{\min} \: \balpha^{*\intercal} \bC \by$, i.e. it is an optimal solution of the multiobjective FOP with weight vector $\balpha^*$. Since $\mathcal{S}^k \subseteq \mathcal{S}$, we then have,
\begin{align*}
    \xi^k_{\text{ABS}}(\balpha^*) &= \balpha^{*\intercal} \bC \hat{\by} - \underset{\by \in \mathcal{S}^k}{\min} \: \balpha^{*\intercal} \bC \by \\
    \xi^k_{\text{ABS}}(\balpha^*) &= \balpha^{*\intercal} \bC \hat{\by} - \balpha^{*\intercal}(\lambda_0 \bC \by^{(0)} + \lambda_1 \bC \by^{(1)} + \dots + \lambda_k \bC \by^{(k)}) \\
    \xi^k_{\text{ABS}}(\balpha^*) &= \xi_{\text{ABS}}(\balpha^*)
\end{align*}
\textit{Proof of (ii)}: Let $I$ denote the the set of indices $i \in 0 \dots k$ for which $\lambda_i > 0$. As stated above, this implies that each $\by^{(i)} \in \mathcal{S}^k, \: i \in I$ is a minimizer of $\underset{\by \in \mathcal{S}^k}{\min} \: \balpha^{*\intercal} \bC \by$. For any $\balpha \in \mathcal{A} \setminus \balpha^*$, we may express $\balpha$ as $\balpha^* + \mathbf{\epsilon}$ for some $\mathbf{\epsilon} \in \mathbb{R}^{|\mathcal{K}|}, \: \mathbf{\epsilon}^\intercal \mathbf{1} = 0, \mathbf{\epsilon} \neq \mathbf{0}$. As such we have,
\begin{align*}
    \xi^k_{\text{ABS}}(\balpha) &= \balpha^{\intercal} \bC \hat{\by} - \underset{\by \in \mathcal{S}^k}{\min} \: \balpha^{\intercal} \bC \by \\
    \xi^k_{\text{ABS}}(\balpha) &= (\balpha^* + \mathbf{\epsilon})^\intercal \bC \hat{\by} - \underset{\by \in \mathcal{S}^k}{\min} \: (\balpha^* + \mathbf{\epsilon})^\intercal \bC \by
\end{align*}
Let $\overline\by$ denote $\underset{\by \in \mathcal{S}^k, i \in I}{\min} \: (\mathbf{\epsilon})^\intercal \bC \by$. Since $\bC \hat{\by} - \xi_{\text{ABS}}(\balpha^*)\mathbf{1}$ lies in the \emph{interior} of a facet of conv($B^k$), then it must be true that there is a set of $|\mathcal{K}|$ linearly independent elements of $\bC \by^{(i)} \in B^k$ such that $i \in I$. Let us denote this set $\overline B$. We note that conv($\overline B$) is a subset of a facet of conv($B^k$), and that $\bC \hat{\by} - \xi_{\text{ABS}}(\balpha^*)\mathbf{1}$ is contained in the interior of this subset. The elements of $\overline B$ define a set of dimension $|\mathcal{K}| - 1$ contained in the plane $\balpha^{*\intercal}\bC\by = (\underset{\by \in \mathcal{S}}{\min} \: \balpha^{*\intercal} \bC \by)$. Since the elements of $\overline B$ are linearly independent and $\balpha^*$ and $\mathbf{\epsilon}$ are not parallel vectors, it must be that there is at least one element $\bC \tilde\by \in \overline B$ such that $\epsilon^\intercal \bC \tilde\by > \epsilon^\intercal \bC \overline\by$. Then it must also be true that $\epsilon^\intercal \bC \overline\by < \epsilon^\intercal (\underset{i \in I}{\sum} \lambda_i \bC \by^{(i)}) = \epsilon^\intercal (\bC \hat{\by} - \xi_{\text{ABS}}(\balpha^*)\mathbf{1})$. With this, we then have that
\begin{align*}
    \xi^k_{\text{ABS}}(\balpha) &\geq \balpha^{*\intercal} \bC \hat{\by} - \balpha^{*\intercal} \bC \overline\by + \mathbf{\epsilon}^\intercal \bC \hat{\by} - \mathbf{\epsilon}^\intercal \bC \overline\by \\
    \xi^k_{\text{ABS}}(\balpha) &\geq \xi^k_{\text{ABS}}(\balpha^*) + \mathbf{\epsilon}^\intercal \bC \hat{\by} - \mathbf{\epsilon}^\intercal \bC \overline\by \\
    \xi^k_{\text{ABS}}(\balpha) &> \xi^k_{\text{ABS}}(\balpha^*) + \mathbf{\epsilon}^\intercal \bC \hat{\by} - \mathbf{\epsilon}^\intercal (\bC \hat{\by} - \xi_{\text{ABS}}(\balpha^*)\mathbf{1})\\
    \xi^k_{\text{ABS}}(\balpha) &> \xi^k_{\text{ABS}}(\balpha^*) + (\mathbf{\epsilon}^\intercal \bC \hat{\by} - \mathbf{\epsilon}^\intercal \bC \hat{\by}) + \xi_{\text{ABS}}(\balpha^*)\mathbf{\epsilon}^\intercal\mathbf{1} \\
    \xi^k_{\text{ABS}}(\balpha) &> \xi^k_{\text{ABS}}(\balpha^*) + \xi_{\text{ABS}}(\balpha^*)*0 \\
    \xi^k_{\text{ABS}}(\balpha) &> \xi^k_{\text{ABS}}(\balpha^*)
\end{align*}
\halmos
\endproof

\subsection{Proof of proposition~\ref{prop:ggtermination}} \label{ECproofGGterminate}

\proof{Proof of proposition~\ref{prop:ggtermination}:} In order to show that the algorithm terminates, we first demonstrate that the algorithm will always reach a termination check criterion in a finite number of iterations, until one such check succeeds. We first show that for any algorithm step $k$, either $\nabla \xi_\text{ABS}(\balpha^{(k)}) = 0$ and the termination check will occur, or there exists a step $k+n$ such that $\bm{\alpha}^{(k+n) \intercal} C \hat{\by} - \bm{\alpha}^{(k+n) \intercal} C \by^{(k+n)} > \xi_\text{ABS} + (\balpha^{(k+n)} - \balpha^{(k+n-1)})^\intercal (C (\hat{\by} - \by^{(k+n-1)}))$ and $\by^{(k+n)} \in \mathcal{S}^{(k+n-1)}$. At a given step $k$, $\balpha^{(k)}$ is found on some facet of the gap function $\xi_\text{ABS}$ defined by $\by^{(k)}$. If $\nabla \xi_\text{ABS}(\balpha^{(k)}) \neq 0$, then the descent method will proceed to take identical steps until some $\balpha^{(k+n)}$ either leaves said facet, in which case $\by^{(k+n)} \neq \by^{(k+n-1)}$, or arrives at the boundary of $\mathcal{A}$ and $\balpha^{(k+n)} = \balpha^{(k+n-1)}$, at which point a termination check will occur. If $\balpha^{(k+n)}$ leaves said facet, then $\by^{(k+n)} \neq \by^{(k+n-1)}$, and since the newfound facet has a different subgradient, the gap function will not decrease by the maximum amount, and thus $\bm{\alpha}^{(k+n) \intercal} C \hat{\by} - \bm{\alpha}^{(k+n) \intercal} C \by^{(k+n)} > \xi_\text{ABS} + (\balpha^{(k+n)} - \balpha^{(k+n-1)})^\intercal (C (\hat{\by} - \by^{(k+n-1)}))$. Since conv$(\mathcal{S})$ has a finite number of extreme points, there exists some step $k+n$ where $\bm{\alpha}^{(k+n) \intercal} C \hat{\by} - \bm{\alpha}^{(k+n) \intercal} C \by^{(k+n)} > \xi_\text{ABS} + (\balpha^{(k+n)} - \balpha^{(k+n-1)})^\intercal (C (\hat{\by} - \by^{(k+n-1)}))$ and $\by^{(k+n)} \in \mathcal{S}^{(k+n-1)}$. Thus, Algorithm~\ref{ggproj} will always reach a termination check in a finite number of steps, and will continue to do so until one succeeds and terminates the algorithm.

Since each failed termination check must necessarily find a vertex of conv($\mathcal{B}$) that is not in $B^k$ and then adds it to $B^{k+1}$, and conv($\mathcal{B}$) has a finite number of vertices, there must be a finite number of failed termination checks until one succeeds, and by Proposition~\ref{lemma:MPsolve}, yields the correct solution. Thus, Algorithm~\ref{ggproj} terminates in a finite number of steps.
\halmos
\endproof

\subsection{Proof of theorem~\ref{thm:ensemble_converge}} \label{ECproofensembleconverge}

\proof{Proof of Theorem~\ref{thm:ensemble_converge}:}
        Let $\mathcal{B}_i = \{\bC\by, \: \forall \by \in \mathcal{S}_i\}$. Since $\xi_\text{ENS}$ lower bounds $\xi_\text{ABS}$, we have by Proposition~\ref{lemma:MPsolve} that if $\mathcal{B}_1 \cup \mathcal{B}_2 \cup \dots \cup \mathcal{B}_n$ contains a set $\bar{B}$ of values of $\bC\by$ whose convex hull contains $\bC\hat{\by} - \xi_\text{ABS}(\balpha^*)$, then $\xi_\text{ENS}$ and $\xi_\text{ABS}$ will have the same minimum value and minimizer. We also know from Proposition~\ref{lemma:MPsolve} that there exists some set $\bar B$ with a total of $|\mathcal{K}|$ elements. Each $\bC\by \in \bar{B}$ is the objective vector of a given partitioning of $G$ into $k$ districts. We note that if the $i^{\text{th}}$ coarsening of a graph $G = \{V, E\}$ does not contract any edges in the edge cut of said districting, then there is a feasible solution to FOP$_i$ that corresponds to the same districting. Since any such districting is a subgraph of $G$ with $L$ disjoint connected components that covers $V$, there are at least $|V|-L$ edges in $E$ do not traverse the edge cut. Further any coarsening of $G$ down to $v$ vertices must contract exactly $|V|-v$ edges. Thus, if $v \geq L$ there exists at least one way of coarsening $G$ such that no edges in the edge cut are contracted. Assuming that our ensemble generation method has a non-zero probability of generating any possible coarsening of $G$ down to $v$ vertices for any coarsening, then for any $\bC\by \in \bar{B}$ and any coarsening $i$ in our ensemble, Prob($\bC\by \not\in \mathcal{B}_i) < 1$. If each coarsening in our ensemble is independently sampled, then as $n \rightarrow \infty$:
        $$\text{Prob}(\bC\by \not\in \mathcal{B}_1 \cup \mathcal{B}_2 \cup \dots \cup \mathcal{B}_n) = \underset{i \in 1 \dots n}{\prod} \text{Prob}(\bC\by \not\in \mathcal{B}_i) \rightarrow 0$$
        
        Since Prob($\bar{B} \not\subseteq \mathcal{B}_1 \cup \mathcal{B}_2 \cup \dots \cup \mathcal{B}_n) = \text{Prob}(\underset{\bb \in \bar{B}}{\bigcup} (\bb \not\in \mathcal{B}_1 \cup \mathcal{B}_2 \cup \dots \cup \mathcal{B}_n))$, we also have that as $n \rightarrow \infty$, Prob($\bar{B} \not\subseteq \mathcal{B}_1 \cup \mathcal{B}_2 \cup \dots \cup \mathcal{B}_n) \rightarrow 0$.
    $\halmos$

\subsection{Proof of proposition~\ref{prop:multipoint_gradient}} \label{ECproofmultipointgradient}

\proof{Proof of Proposition~\ref{prop:multipoint_gradient}:}
    If we can find a subtangent plane of $\xi_\text{ENS}(\balpha)$ at $\balpha^{(k)}$, then the gradient of said plane is a subgradient of $\xi_\text{ENS}(\balpha)$. Let us define our subtangent plane at $\balpha^{(k)}$ as $P(\balpha) =  \balpha^\intercal \bC \hat{\by} - \balpha^\intercal \bC (\underset{\by \in \by^{(k)}_1 \dots \by^{(k)}_n}{\argmin} \: \balpha^{(k)\intercal} \bC \by)$. Then, $\nabla P(\balpha) = \bC (\hat\by - \underset{\by \in \by^{(k)}_1 \dots \by^{(k)}_n}{\argmin} (\balpha^{(k)\intercal}\bC\by))$. To show that $P(\balpha)$ is subtangent at $\balpha^{(k)}$, we show that (i) $P(\balpha^{(k)}) = \xi_\text{ENS}(\balpha^{(k)})$, and (ii) $P(\balpha) \leq \xi_\text{ENS}(\balpha)$.

    \textit{Proof of (i)}: we note by the formulation of MultiGIO$_\text{ABS}$MinMax that $\xi_\text{ENS}(\balpha) = \underset{i \in 1 \dots n}{\max}\xi_i(\balpha)$ where $\xi_i(\balpha)$ denotes the gap function generated by the inverse formulation GIO$_\text{ABS}$ applied to FOP$_i$. Then we have that:
\begin{align*}
    \xi_\text{ENS}(\balpha) &= \underset{i \in 1 \dots n}{\max}\xi_i(\balpha) \\
    \xi_\text{ENS}(\balpha) &= \underset{i \in 1 \dots n}{\max}(\balpha^\intercal \bC \hat{\by} - \underset{\by \in \mathcal{S}_i}{\min}\balpha^\intercal \bC \by)
\end{align*}
By Proposition~\ref{thm:gap_func}, we have that:
\begin{align*}
    \xi_\text{ENS}(\balpha^{(k)}) &= \underset{i \in 1 \dots n}{\max}(\balpha^{(k)\intercal} \bC \hat{\by} - \balpha^{(k)\intercal} \bC \by^{(k)}_i) \\
    \xi_\text{ENS}(\balpha^{(k)}) &= \balpha^{(k)\intercal} \bC \hat{\by} - \balpha^{(k)\intercal}\bC(\underset{\by \in \by^{(k)}_1 \dots \by^{(k)}_n}{\argmin} \: \balpha^{(k)\intercal} \bC \by) \\
    \xi_\text{ENS}(\balpha^{(k)}) &= P(\balpha^{(k)})
\end{align*}
    \textit{Proof of (ii)}:
\begin{align*}
    \xi_\text{ENS}(\balpha) &= \underset{i \in 1 \dots n}{\max}(\balpha^\intercal \bC \hat{\by} - \underset{\by \in \mathcal{S}_i}{\min}\balpha^\intercal \bC \by) \\
    \xi_\text{ENS}(\balpha) &= \balpha^\intercal \bC \hat{\by} - \underset{\by \in \mathcal{S}_1 \dots \mathcal{S}_n}{\min}\balpha^\intercal \bC \by \\
    \xi_\text{ENS}(\balpha) &\geq \balpha^{\intercal} \bC \hat{\by} - \balpha^{\intercal}\bC(\underset{\by \in \by^{(k)}_1 \dots \by^{(k)}_n}{\argmin} \: \balpha^{(k)\intercal} \bC \by) \\
    \xi_\text{ENS}(\balpha) &\geq P(\balpha)
\end{align*}
\halmos

\section{Relative sub-optimality loss function}\label{sec:EC_MTrelgap}

The relative sub-optimality loss function, $\ell_{\text{rel}}(\hat\by, \mathcal{S}, \balpha) = \frac{\balpha^\intercal \bC \hat\by}{\underset{\by \in \mathcal{S}}{\min}\balpha^\intercal \bC \by}$ measures the quotient of the input solution objective value and the optimal objective value. We note that this loss function is well defined only if $\balpha \in \mathbb{R}_+^{|\mathcal{K}|} \setminus \{\bold{0}\}$ and $\bC \by > \bold{0}, \forall \by \in \mathcal{S}$. The corresponding data-driven inverse optimization problem with relative sub-optimality as a loss function can equivalently be written as:
\begin{align*}
    \minimize_{\bm{\alpha},\; \xi_{\text{rel}}} \quad & \xi_{\text{rel}} & \tag{GIO$_{\text{rel}}$} \\
    \text{subject to } \quad&\bm{\alpha}^\intercal \bC \hat\by \leq (\bm{\alpha}^\intercal \bC \by) \xi_{\text{rel}}, 
    \quad \forall \: \by \in \mathcal{S},\\
    &\balpha \in \mathcal{A},\\
    &\bm{\alpha} \geq 0.
\end{align*}

Note that GIO$_{\text{rel}}$ is not a linear program because two continuous decision variables are multiplied ($\bm{\alpha}$ and $\xi_{\text{rel}}$) in the first set of constraints. \cite{chan2014generalized} solve a similar formulation by varying $\xi_{\text{rel}}$ with a univariate search technique until the smallest value for $\xi$ is found such that GIO$_{\text{rel}}$ is feasible. However in the case that $\mathcal{A}$ is defined as the unit simplex, using a method similar to \cite{chan2019inverse}, the exact minimum relative gap and corresponding objective weighting can be derived quickly from the solution of a single linear program. 

\begin{theorem}\label{thm1}
$\{\bm{\bar{\alpha}}, \bar{\xi}_{\text{rel}}\}$ is an optimal solution to \emph{(GIO$_{\text{rel}}$)} if $\bm{\bar{\alpha}} = \frac{\bm{\hat{\alpha}}}{||\bm{\hat{\alpha}}||_1}$ and $\bar{\xi}_{\text{rel}} = \bm{\hat{\alpha}}^\intercal \bC \hat\by$, where $\bm{\hat{\alpha}}$ is an optimal solution to the following linear program \emph{(GIO$_{\text{rel}}^*$)}:
\begin{align*}
    \minimize_{\bm{\alpha}} \quad & \bm{\alpha}^\intercal \bC \hat\by & \tag{GIO$_{\text{rel}}^*$} \\
    \text{subject to } \quad &\bm{\alpha}^\intercal \bC \by \geq 1, \quad \forall \: \by \in \mathcal{S},\\
    &\bm{\alpha} \geq 0.
\end{align*}
\end{theorem}

\proof{Proof of Theorem~\ref{thm1}:} 
First, we show that (i) $\{\bm{\bar{\alpha}}, \bar{\xi}_{\text{rel}}\}$ is a feasible solution to (GIO$_{\text{rel}}$). Next, we show that (ii) if $\{\bm{\bar{\alpha}}, \bar{\xi}_{\text{rel}}\}$ is not optimal for (GIO$_{\text{rel}}$), then $\bm{\hat{\alpha}}$ is not optimal for for (GIO$_{\text{rel}}^*$), thus proving the claim by contrapositive.

(i) We observe that $||\bm{\bar{\alpha}}||_1 = 1$ and $\bm{\bar{\alpha}} \geq \bzero$ by the construction of $\bm{\bar{\alpha}}$ and $\bar{\xi}_{\text{rel}} = \bm{\hat{\alpha}}^\intercal C \hat\by \geq 0$ because $C \hat\by \in \mathbb{R}^{|\mathcal{K}|}_+$. The remaining constraints are satisfied because $\forall \by \in \mathcal{S}$:

\begin{equation*}
    \bm{\bar{\alpha}}^\intercal C \hat\by = \frac{\bm{\hat{\alpha}}^\intercal C \hat\by}{||\bm{\hat{\alpha}}||_1} = \frac{\bar{\xi}_{\text{rel}}}{||\bm{\hat{\alpha}}||_1} \leq (\bm{\hat{\alpha}}^\intercal C \by) \frac{\bar{\xi}_{\text{rel}}}{||\bm{\hat{\alpha}}||_1} = \frac{\bm{\hat{\alpha}}^\intercal C \by}{||\bm{\hat{\alpha}}||_1} \bar{\xi}_{\text{rel}} = (\bm{\bar{\alpha}}^\intercal C \by) \bar{\xi}_{\text{rel}}.
\end{equation*}

(ii) Suppose there exists a feasible solution $\{\bm{\alpha^*}, \xi_{\text{rel}}^*\}$ to (GIO$_{\text{rel}}$) such that $\xi_{\text{rel}}^* < \bar{\xi}_{\text{rel}}$, i.e. $\{\bm{\bar{\alpha}}, \bar{\xi}_{\text{rel}}\}$ is not an optimal solution. Let $\bm{\tilde{\alpha}} = \frac{\xi_{\text{rel}}^* \bm{\alpha^*}}{\bm{\alpha}^{* \intercal} C \hat\by}$. Then, $\bm{\tilde{\alpha}}$ is a feasible solution to (GIO$_{\text{rel}}^*$) because $\bm{\tilde{\alpha}} \geq \bm{0}$ and $\forall \by \in \mathcal{S}$: 
$$\bm{\tilde{\alpha}}^\intercal C \by = \frac{\xi_{\text{rel}}^*}{\bm{\alpha}^{* \intercal} C \hat\by} (\bm{\alpha}^{* \intercal} C \by) \geq \frac{\xi_{\text{rel}}^*}{\bm{\alpha}^{* \intercal} C \hat\by} \frac{\bm{\alpha}^{* \intercal} C \hat\by}{\xi_{\text{rel}}^*} = 1.$$ 

Since $\bm{\tilde{\alpha}}^\intercal C \hat\by = \xi_{\text{rel}}^* < \bar{\xi}_{\text{rel}} = \bm{\hat{\alpha}}^\intercal C \hat\by$, $\bm{\hat{\alpha}}$ is not an optimal solution to (GIO$_{\text{rel}}^*$), which contradicts the definition of $\bm{\hat{\alpha}}$. Therefore, such an $\{\bm{\alpha^*}, \xi_{\text{rel}}^*\}$ cannot exist, so $\{\bm{\bar{\alpha}}, \bar{\xi}_{\text{rel}}\}$ must be an optimal solution for (GIO$_{\text{rel}}$). 
\halmos
\endproof

Theorem~\ref{thm1} shows that we can obtain the optimal solution to GIO$_{\text{rel}}$ (which is bilinear) by solving GIO$_{\text{rel}}^*$ (which is a linear program).

As such, our modification to the algorithm proposed by \cite{moghaddass2020inverse} minimizes the relative sub-optimality loss function, rather than absolute sub-optimality. The algorithm structure is shown in Algorithm~\ref{alg_MTrelgap}. The main distinction is that our algorithm solves GIO$_{\text{rel}}^*$ and uses Theorem~\ref{thm1} to obtain the solution to GIO$_{\text{rel}}$. 

\begin{theorem}\label{thm2}
    Algorithm~\ref{alg_MTrelgap} terminates finitely with an optimal solution to GIO$_{\text{rel}}$.
\end{theorem}

\proof{Proof of Theorem~\ref{thm2}:} 
The proof follows directly from the proof of termination and correctness in \cite{moghaddass2020inverse}, together with the proof of Theorem~\ref{thm1}.
\halmos
\endproof

Theorem~\ref{thm2} demonstrates that Algorithm~\ref{alg_MTrelgap} produces the optimal solution to GIO$_{\text{rel}}$ when the complete forward problem feasible region is known, in at most as many steps as there are extreme points to the forward multi-objective feasible region.

\begin{algorithm}
 \SetKwInOut{Input}{input}\SetKwInOut{Output}{output}
 \Input{$C, \: \bm{y}^0$, FP}
 \Output{$\bm{\alpha}^\text{best}, \: \xi^\text{best}_\text{rel}$}
 $\tilde{\mathcal{S}} \leftarrow \emptyset, \: \bm{\hat{\alpha}} \leftarrow \text{GenInvOp}_{\text{rel}}^*(\bm{y}^0, \tilde{\mathcal{S}})$\;
 $\bm{\alpha} \leftarrow \frac{\bm{\hat{\alpha}}}{||\bm{\hat{\alpha}}||_1}, \: \bm{y} \leftarrow \text{FP}(\bm{\alpha})$\;
 $\xi^\text{best}_\text{rel} \leftarrow \bm{\alpha}^\intercal C \bm{y}^0 - \bm{\alpha}^\intercal C \bm{y}, \: \bm{\alpha}^\text{best} \leftarrow \bm{\alpha}$\;
 \While{$\bm{\alpha}^\intercal C \bm{y}^0 > \bm{\alpha}^\intercal C \bm{y}$}{
  $\tilde{\mathcal{S}} \leftarrow \tilde{\mathcal{S}} \cup \bm{y}$\;
  $\bm{\hat{\alpha}} \leftarrow \text{GenInvOp}_{\text{rel}}^*(\bm{y}^0, \tilde{\mathcal{S}})$\;
  $\bm{\alpha} \leftarrow \frac{\bm{\hat{\alpha}}}{||\bm{\hat{\alpha}}||_1}$\;
  $\bm{y} \leftarrow \text{FP}(\bm{\alpha})$\;
  \If{$\bm{\alpha}^\intercal C \bm{y}^0 - \bm{\alpha}^\intercal C \bm{y} = \xi^\text{best}_\text{rel}$}{
   \eIf{$\bm{\alpha} = \bm{\alpha}^\text{best}$}{
    \textbf{stop}
    }{
    $\bm{\alpha}^\text{best} \leftarrow \bm{\alpha}$\;
    }
  \If{$\bm{\alpha}^\intercal C \bm{y}^0 - \bm{\alpha}^\intercal C \bm{y} < \xi^\text{best}_\text{abs}$}{
   $\xi^\text{best}_\text{abs} \leftarrow \bm{\alpha}^\intercal C \bm{y}^0 - \bm{\alpha}^\intercal C \bm{y}$\;
   $\bm{\alpha}^\text{best} \leftarrow \bm{\alpha}$\;
  }
 }
}    
 \caption{Data-driven relative gap cutting-plane method.}
 \label{alg_MTrelgap}
\end{algorithm}

\section{Adapting gap-gradient methods to the relative sub-optimality loss function} \label{sec:relgap_adaptation}
In the case of the relative optimality gap, we apply gap-gradient methods to an inverse formulation based on the linearization $GIO_\text{REL}^*$. Thus, for relative gap problems, we can formulate the corresponding gap function as $\xi_\text{REL}(\balpha) = \balpha^\intercal\bC\hat\by$ defined over the gap function domain $\mathcal{A} = \{\balpha \in \mathbb{R}^\mathcal{K} \bigm\lvert \balpha^\intercal C \by \geq 1, \balpha \geq 0, \forall \by \in \mathcal{S}\}$

Similar to the absolute gap case, we note that the feasible region of the master problem GIO$_\text{rel}^*$ in the case where all extreme points are known is the epigraph of the relative gap function, and we obtain the following result.  
\begin{proposition}\label{thm:gap_rel}
    The relative gap function $\xi_\text{REL}(\balpha)$ is a convex function.
\end{proposition}
\proof{Proof of Theorem~\ref{thm:gap_rel}:}
The gap function is a linear function, and as such is convex if its domain is convex. The domain of the gap function is the union of half-planes, and as such, is convex.
\halmos
\endproof

\begin{figure}
    \centering
    \begin{tikzpicture}[scale=1]
    \draw [->] (0, 0) -- (0, 5);
    \draw [->] (0, 0) -- (5, 0);
    \draw [pink] [fill=pink] (1, 5) -- (1, 3.4) -- (1.8, 1.7) -- (3.6, 1) -- (5, 1) -- (5,5);
    \draw [red] (1, 5) -- (1, 3.4) -- (1.8, 1.7) -- (3.6, 1) -- (5, 1);
    \draw [cyan] [-stealth] (2.7, 1.8) -- (3.1, 2.85) node [right] {$\balpha^*$};
    \filldraw [blue] (2.7, 1.8) circle (1pt) node [right] {$\hat\by$};
    \node at (5, -0.4) {$c_1{}^\intercal y$};
    \node at (-0.6, 5) {$c_2{}^\intercal y$};
    \node at (2.5, -1) {(a)};
    \end{tikzpicture}
    \hspace{5pt}
    \begin{tikzpicture}[scale=1]
    \draw [->] (0, 0) -- (0, 5);
    \draw [->] (0, 0) -- (5, 0);
    \draw [pink] [fill=pink] (0, 5) -- (0.8, 2.1) -- (1.9, 0.9) -- (5, 0) -- (5,5);
    \draw [red] (0, 5) node [above right] {domain of $\xi_{\text{REL}}$} -- (0.8, 2.1) -- (1.9, 0.9) -- (5, 0);
    \draw [blue] [-stealth] (0, 0) -- (2.7, 1.8) node [right] {$\bC\hat\by$};
    \filldraw [cyan] (0.8, 2.1) circle (1pt) node [above right] {$\balpha^*$};
    \node at (5, -0.4) {$\balpha_1$};
    \node at (-0.4, 5) {$\balpha_2$};
    \node at (2.5, -1) {(c)};
    \end{tikzpicture}
    
    \caption{An example inverse optimization problem with two sub-objectives. (a) The FOP objective feasible space and inverse input $\hat\by$, and (b) the corresponding relative gap function projected onto two dimensions, and the gradient of the gap function $\nabla\xi_{\text{REL}}(\balpha)=\bC\hat\by$.}
    \label{fig:gap_space_example_rel}
\end{figure}
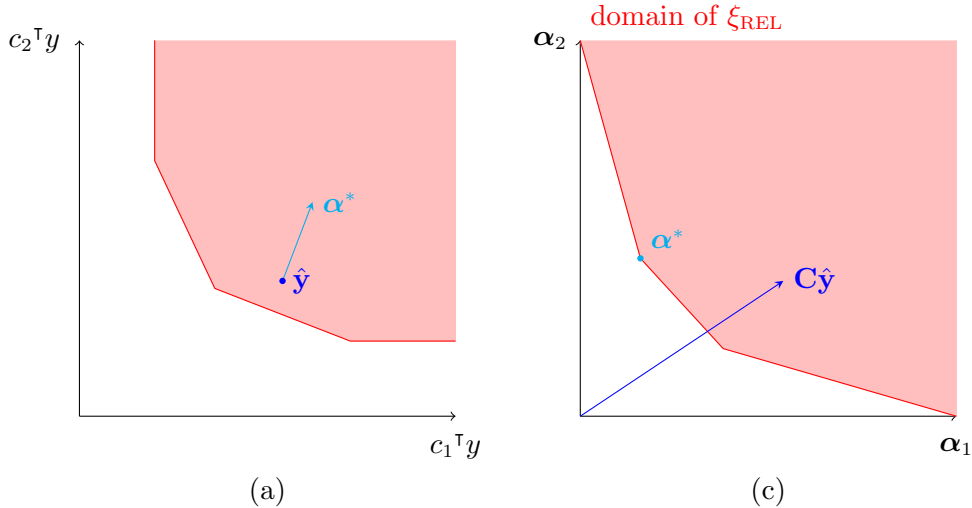

Figure~\ref{fig:gap_space_example_rel}(b) displays the domain of $\xi_{\text{REL}}$ for an example inverse optimization problem with two objectives, and the vector $\bC \hat\by$ which is equal to the gradient of the gap function $\nabla \xi_{\text{REL}}(\balpha)$. Minimizing the gap function yields the minimizer $\balpha^*$.

\begin{remark}
Under the formulation of the relative gap space and relative gap function, each FOP extreme point $\by^{(k)}$ does not yield a tangent plane, but a facet of the boundary of the domain of the relative gap function. For the relative gap function the gradient is in fact uniform ($\nabla g_{\text{REL}}(\balpha^{(k)}) = \bC \hat\by, \forall \balpha^{(k)}$ in the interior of the relative gap function's domain), but each facet gained from an extreme point can be used to find the projection of the gradient onto the boundary of the relative gap function domain at that point. Specifically, for a facet-defining boundary plane $\balpha^\intercal C \by^{(k)} = 1$ discovered by solving the forward problem at objective weighting $\balpha^{(k)}$, the gradient projected onto this boundary is $(\nabla g_{\text{REL}}(\balpha^{(k)}))_\text{proj} = C \hat\by - (\frac{(C \hat\by)^\intercal C \by^{(k)}}{||C \by^{(k)}||^2_2}) C \by^{(k)}$. Similar to the absolute sub-optimality loss function, this can be used to minimize the convex relative gap function in a bounded space.
\end{remark}

A gap-gradient solution method that uses this subgradient could proceed at each iteration by taking a step in the direction of the projected gradient. While it is possible that the $\balpha^{(k)}$ generated by this step lies outside of the domain $\mathcal{A}$ due to incomplete knowledge of its boundary, we note that so long as $\balpha^{(k)}$ lies within the positive orthant, solving FOP$(\balpha^{(k)})$ will yield the same solution $\by^{(k)}$ as an FOP with some scalar-multiplied objective weight $m\balpha^{(k)}, \: m > 1$ such that $m\balpha^{(k)} \in \mathcal{A}$. Once $\by^{(k)}$ is known then the fact-defining plane $\balpha^\intercal \bC \by^{(k)} = 1$ can be used to project $\balpha^{(k)}$ onto the domain boundary by way of the projection $m = \frac{1}{\balpha^{(k)^\intercal} \bC \by^{(k)}}$. One can then proceed in the gap-gradient algorithm by taking a step from $m\balpha^{(k)}$ in the direction of the next projected (negative) gradient, $-\bC\hat\by + (\frac{(C \hat\by)^\intercal C \by^{(k)}}{||C \by^{(k)}||^2_2}) C \by^{(k)}$.

\section{Frank-Wolfe gap gradient methods with FOP early stopping} \label{FWpartialSolve}

\begin{figure}
    \centering
    \includegraphics[scale = 0.6]{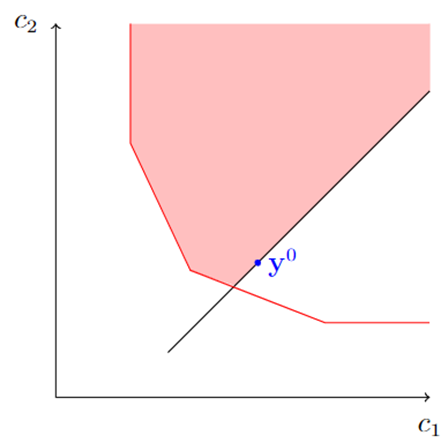}
    \includegraphics[scale = 0.6]{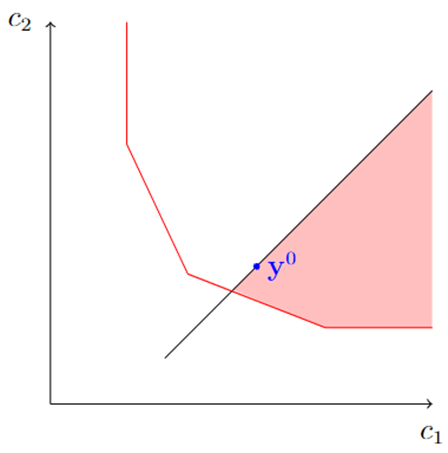}
    \caption{The feasible regions of two concurrently partially solved variations of the FOP shown in Figure \ref{fig:gap_space_example} (a), used in conjunction with the proposed Frank-Wolfe method for potentially faster descent steps.}
    \label{fig:FWshortcut}
\end{figure}

To operationalize this method of applying early stopping to solving the FOP when sufficient information is known, we can use the following approach. First, we create $|\mathcal{K}|$ different variations of the forward problem, each with an additional set of constraints $\bC_k (\by - \hat{\by}) \leq \bC_l (\by - \hat{\by}), \: \forall l \in 1 \dots \mathcal{K} \setminus k$, which enforces that subobjective component $k$ of any feasible solution minus that of the input solution must be less than or equal to that of every other subobjective component. Figure \ref{fig:FWshortcut} illustrates the FOP variations created by these constraints for the example two-objective inverse optimization problem shown in Figure \ref{fig:gap_space_example} (a). Next, we solve each of the FOP$_k$, $k\in\mathcal{K}$ in parallel. Once we reach a point where the lower objective bound of a given FOP$_j$, $j\in\mathcal{K}$, is greater than the lowest upper bound of any other running FOP$_k$, we can terminate the solving process of FOP$_j$. When only one such FOP$_k$ remains, we terminate because it is unnecessary to solve to optimality because we know that whatever the optimal solution is, the gradient will have subobjective $k$ as its greatest component, which is sufficient information to execute the next step in algorithm (the next descent step is in the direction towards $\bm{e}_k$). This concurrent solving subalgorithm is detailed in the Appendix as Algorithm~\ref{ggfwpartialsub}.

Using this method, we can also, at the point of terminating the final running subproblem, return the incumbent solution at the time of termination. This returns a feasible forward solution that produces a lower-bounding plane of the gap function. If it is an interior point, then it will not be tangent to the gap function, but strictly below. However, as our Frank-Wolfe method approaches $\balpha^*$, the concurrent process will be terminating in smaller and smaller FOP optimality gaps, so the returned incumbents will be progressively closer to the boundary of conv($\mathcal{B}$). Thus, once we are within a close neighborhood of the gap function minimizer, we may collect a set of proper tangent planes (at this point our method is likely to be returning vertices of conv($\mathcal{B}$)), and we can use these gathered tangent planes in a termination check that is similar to the termination methods used for our other proposed methods. The methods discussed in this section currently remain unimplemented because they will likely involve constructing an MILP solver from the ground up.

\section{Formulation.}\label{app:FOP-formulation}

The overall mixed-integer linear formulation for the FOP is given by:
\allowdisplaybreaks
\begin{subequations} \label{FP}
\begin{align}
    &\minimize_{x, z^D, v^D, f, \rho, \sigma_A, \phi_{EG}} &\alpha_1 \rho + \alpha_2 \sigma_A + \alpha_3 \phi_{EG}& & &\\
    &\text{subject to } &\sum_{i \in V} x_{ii} &= L, & &\\
    &&\sum_{i \in V} x_{ij} &= 1, &\forall \: j &\in V, \\
    &&x_{ij} &\leq x_{ii}, &\forall \: i, j &\in V, \\
    &&x_{ij} + \sum_{v \in N(j)} (f_{ijv} - f_{ivj}) &= 0, &\forall i, j &\in V, i \neq j, \\
    &&x_{ii} + \sum_{v \in N(i)} (f_{iiv} - f_{ivi}) - \sum_{v \in N(i)} x_{iv} &= 0, &\forall i &\in V, \\
    &&|V|x_{ij} - \sum_{v \in N(j)} f_{ivj} &\geq 0, &\forall i, j &\in V, \\
    &&(x_{ii} - \rho) \overline{P} &\leq \sum_{j \in V} p_j x_{ij}, &\forall i &\in V, \\
    &&(x_{ii} + \rho) \overline{P} &\geq \sum_{j \in V} p_j x_{ij}, &\forall i &\in V, \\
    &&\sigma_A &= \frac{\sum_{i, j \in V}d_{ij} a_j x_{ij}}{M}, \\
    &&\overline{P} \leq \sum_{j \in V}(p^D_j - p^R_j)x_{ij}, - \overline{P}z^D_i &\leq 0, &\forall i &\in V, \\
    &&v^D_{ij} &\leq x_{ij}, &\forall i, j &\in V, \\
    &&v^D_{ij} &\leq z^D_{i}, &\forall i, j &\in V, \\
    &&v^D_{ij} &\geq x_{ij} + z^D_{i} - 1, &\forall i, j &\in V, \\
    &&\sum_{j \in V} (\frac{3p^D_j - p^R_j}{2}) x_{ij} - \sum_{j \in V} (p^D_j + p^R_j) v^D_{ij} &= w_i, &\forall i &\in V, \\
    &&\phi_{EG} &\geq \frac{\sum_{i \in V} w_i}{\sum_{i \in V}(p^D_i + p^R_i)}, \\
    &&\phi_{EG} &\geq -\frac{\sum_{i \in V} w_i}{\sum_{i \in V}(p^D_i + p^R_i)}, \\
    && x_{ij}, z^D_i, v^D_{ij} \in \{0,1\}, \: f_{ivj} &\geq 0,  &\forall i, j &\in V, \forall v \in N(j).
\end{align}
\end{subequations}

Constraints (2b) through (2d) enforce that every vertex is assigned to exactly one district out of $L$ total districts. Constraints (2e) through (2g) construct a set of flow networks that maintain contiguity for the created districts. Constraints (2h) and (2i) define the $\rho$ variable used in the objective. Constraint (2j) defines $\sigma_A$ (used in the objective). Finally, Constraints (2k) - (2q) define $\phi_{EG}$, which is used in the objective. 

\section{Experimental details}

\subsection{MIPLIB instances} \label{MIPLIBfops}

The following are labels for the FOPs used in Section~\ref{MIPLIBExperiments}:

\begin{enumerate}
    \item neos-1430701

    \item gsvm2rl3

    \item ran13x13

    \item spd150x300d

    \item supportcase17


    \item ci-24

    \item ic97\_tension

    \item fastxgemm-n2r60t2

    \item timtab1CUTS
    
\end{enumerate}

\subsection{Generating sample states}\label{sec:exp}

To create a sample state $G$ for a specified $|V|$, we randomly sample $|V|$ points uniformly over a unit square, and calculate the Delaunay triangulation to create an adjacency graph of our simulated census blocks. The distance matrix is calculated using the euclidean distances of the sampled points. Note that scaling of the distance matrix does not matter, as the compactness metric of a districting is normalized for the 1-median of the entire state. Land areas for the corresponding simulated census blocks are calculated as the areas of the Voronoi cells produced by the sampled points. At each census block, values for the number of voters for party A, party B, and non-voting people are randomly sampled from integers in the range [10, 100]. For each sampled state, an FOP is generated for the political districting optimization model with $L = 2$ districts.

\section{Algorithm Structures}\label{sec:EC_algstruc}

\begin{algorithm}[H]
 \SetKwInOut{Input}{input}\SetKwInOut{Output}{output}
 \Input{$\bC, \: \hat\by$, FOP}
 \Output{$\bm{\alpha}^\text{best}, \: \xi^\text{best}$}
 $\tilde{\mathcal{S}} \leftarrow \emptyset, \: \bm{\alpha} \leftarrow \text{GIO}_{\text{ABS}}(\hat\by, \tilde{\mathcal{S}}), \: \by \leftarrow \text{FOP}(\bm{\alpha})$\;
 $\xi^\text{best} \leftarrow \bm{\alpha}^\intercal \bC \hat\by - \bm{\alpha}^\intercal \bC \by, \: \bm{\alpha}^\text{best} \leftarrow \bm{\alpha}$\;
 \While{$\bm{\alpha}^\intercal \bC \hat\by > \bm{\alpha}^\intercal C \by$}{
  $\tilde{\mathcal{S}} \leftarrow \tilde{\mathcal{S}} \cup \by$\;
  $\bm{\alpha} \leftarrow \text{GIO}_{\text{ABS}}(\hat\by, \tilde{\mathcal{S}})$\;
  $\by \leftarrow \text{FP}(\bm{\alpha})$\;
  \If{$\bm{\alpha}^\intercal C \hat\by - \bm{\alpha}^\intercal C \by = \xi^\text{best}$}{
   \eIf{$\bm{\alpha} = \bm{\alpha}^\text{best}$}{
    \textbf{stop}
    }{
    $\bm{\alpha}^\text{best} \leftarrow \bm{\alpha}$\;
    }
  \If{$\bm{\alpha}^\intercal C \hat\by - \bm{\alpha}^\intercal C \by < \xi^\text{best}_\text{abs}$}{
   $\xi^\text{best}_\text{abs} \leftarrow \bm{\alpha}^\intercal C \hat\by - \bm{\alpha}^\intercal C \by$\;
   $\bm{\alpha}^\text{best} \leftarrow \bm{\alpha}$\;
  }
 }
}
 \caption{Moghaddass and Terekhov (2020)} \label{moghaddass}
\end{algorithm}

\begin{algorithm}
 \SetKwInOut{Input}{input}\SetKwInOut{Output}{output}
 \Input{$C, \: \hat{\by}$, FOP, $t$, $\beta$}
 \Output{$\bm{\alpha}^\text{best}, \: \xi^\text{final}$}
 $k = 0, \: \mathcal{S}^k \leftarrow \emptyset, \: \bm{\alpha}^{(k)} \leftarrow \frac{1}{\mathcal{K}} \mathbf{1}$\;
 $\by^{(k)} \leftarrow \text{FOP}(\bm{\alpha}^{(k)})$\;
 $\xi_\text{ABS} \leftarrow \bm{\alpha}^{(k) \intercal} C \hat{\by} - \bm{\alpha}^{(k) \intercal} C \by^{(k)}$\;
 \While{True}{
  $k \leftarrow k + 1$\;
  $\mathcal{S}^k \leftarrow \mathcal{S}^{k-1} \cup \by^{(k-1)}$\;
  $\bm{\alpha}^{(k)} \leftarrow \bm{\alpha}^{(k-1)} - t (C (\hat{\by} - \by^{(k-1)})) + \beta(\balpha^{(k-1)} - \balpha^{(k-2)})$\;
  $\bm{\alpha}^{(k)} \leftarrow \text{proj}_{\Delta^\mathcal{K}}(\bm{\alpha}^{(k)})$\;
  $\by^{(k)} \leftarrow \text{FOP}(\bm{\alpha}^{(k)})$\;
  \eIf{$\balpha^{(k)} = \balpha^{(k-1)}$\textbf{or} $\bm{\alpha}^{(k) \intercal} C \hat{\by} - \bm{\alpha}^{(k) \intercal} C \by^{(k)} > \xi_\text{ABS} + (\balpha^{(k)} - \balpha^{(k-1)})^\intercal (C (\hat{\by} - \by^{(k-1)}))$ \textbf{and} $\by^{(k)} \in \mathcal{S}^k$}{
    $\balpha^{\text{final}}, \xi^{\text{final}} \leftarrow \text{GIO}_{\text{ABS}}(\hat{\by}, \mathcal{S}^k)$\;
    $k \leftarrow k + 1$\;
    $\by^{(k)} \leftarrow \text{FOP}(\bm{\alpha}^{\text{final}})$\;
    \eIf{$\balpha^{\text{final} \intercal} C \hat\by - \balpha^{\text{final} \intercal} C \by^{(k)} = \xi^\text{final}$}{
    \textbf{stop}
    }{
    $\balpha^{(k)} \leftarrow \balpha^{\text{final}}$\;
    $\mathcal{S}^k \leftarrow \mathcal{S}^{k-1} \cup \by^{(k-1)}$\;
    $t \leftarrow \frac{t}{2}$\;
    }
  }{
   $\xi_\text{ABS} \leftarrow \bm{\alpha}^{(k) \intercal} C \hat{\by} - \bm{\alpha}^{(k) \intercal} C \by^{(k)}$\;
 }
}    
 \caption{Gap-gradient Projected Gradient Descent with Acceleration.}
 \label{ggheavyball}
\end{algorithm}

\begin{algorithm}
 \SetKwInOut{Input}{input}\SetKwInOut{Output}{output}
 \Input{$C, \: \hat\by$, FOP}
 \Output{$\balpha^\text{final}, \: \xi^\text{final}$}
 $k = 0, \: \mathcal{S}^k \leftarrow \emptyset, \: \balpha^{(k)} \leftarrow \frac{1}{\mathcal{K}} \mathbf{1}$\;
 $\by^{(k)} \leftarrow \text{FOP}(\bm{\alpha}^{(k)})$\;
 $\xi_\text{abs} \leftarrow \bm{\alpha}^{(k) \intercal} C \hat\by - \bm{\alpha}^{(k) \intercal} C \by^{(k)}$\;
 \While{True}{
  $k \leftarrow k + 1$\;
  $\mathcal{S}^k \leftarrow \mathcal{S}^{k-1} \cup \by^{(k-1)}$\;
  $i \leftarrow \argmin_i C_i(\hat\by - \by^{(k-1)})$\;
  $\bm{\alpha}^{(k)} \leftarrow \bm{\alpha}^{(k-1)} - \frac{2}{2+k}(\bm{e}_i - \balpha^{(k-1)})$\;
  $\by^{(k)} \leftarrow \text{FOP}(\bm{\alpha}^{(k)})$\;
  \eIf{$\balpha^{(k)} = \balpha^{(k-1)}$\textbf{or} $\bm{\alpha}^{(k) \intercal} C \hat{\by} - \bm{\alpha}^{(k) \intercal} C \by^{(k)} > \xi_\text{ABS} + (\balpha^{(k)} - \balpha^{(k-1)})^\intercal (C (\hat{\by} - \by^{(k-1)}))$ \textbf{and} $\by^{(k)} \in \mathcal{S}^k$}{
    $\balpha^{\text{final}}, \xi^{\text{final}} \leftarrow \text{GIO}_{\text{ABS}}(\hat{\by}, \mathcal{S}^k)$\;
    $k \leftarrow k + 1$\;
    $\by^{(k)} \leftarrow \text{FOP}(\bm{\alpha}^{\text{final}})$\;
    \eIf{$\balpha^{\text{final} \intercal} C \hat\by - \balpha^{\text{final} \intercal} C \by^{(k)} = \xi^\text{final}$}{
    \textbf{stop}
    }{
    $\balpha^{(k)} \leftarrow \balpha^{\text{final}}$\;
    $\mathcal{S}^k \leftarrow \mathcal{S}^{k-1} \cup \by^{(k-1)}$\;
    }
  }{
   $\xi_\text{ABS} \leftarrow \bm{\alpha}^{(k) \intercal} C \hat{\by} - \bm{\alpha}^{(k) \intercal} C \by^{(k)}$\;
 }
}
\caption{Frank-Wolfe Generalized Inverse Method, Absolute Gap}
\label{ggfw}
\end{algorithm} 

\begin{algorithm}
 \SetKwInOut{Input}{input}\SetKwInOut{Output}{output}
 \Input{FOP, $\bC$, $\bm{\alpha}$, $\bm{y}^0$}
 \Output{$i$, $\bm{y}$}
 \For{$i \in 1 \dots \mathcal{K}$}{ $\text{FOP}_i \leftarrow \text{FOP} \cap \{C_i (\bm{y} - \bm{y^0}) \leq C_j (\bm{y} - \bm{y^0}) \: \big{\vert} \: j \in 1 \dots \mathcal{K} \setminus i\}$\;
 }
 $u \leftarrow \infty$\;
 $p \leftarrow \mathcal{K}$\;
 $\forall i$, concurrently initiate optimizing FOP$_i$ with objective coefficients $\bm{\alpha}$\;
 \While{FOP$_i$ is solving}{
  \If{FOP$_i$.upper\_bound $< u$}{
    $u \leftarrow$ FOP$_i$.upper\_bound\;
  }
  \If{FOP$_i$.lower\_bound $> u$}{
    $p \leftarrow p - 1$\;
    terminate FOP$_i$\;
  }
  \If{$p = 1$}{
    $\bm{y} \leftarrow$ incumbent solution of FOP$_i$\;
    terminate FOP$_i$\;
    \Return $i$, $\bm{y}$
  }
}    
 \caption{Frank-Wolfe Incomplete FOP Solve Subalgorithm, Concurrent MIP Approach, Absolute Gap}
 \label{ggfwpartialsub}
\end{algorithm} 

\begin{algorithm}
 \SetKwInOut{Input}{input}\SetKwInOut{Output}{output}
 \Input{$\bC, \: \hat\by$, FOP}
 \Output{$\bm{\alpha}^\text{final}, \: \xi^\text{final}$}
 $k = 0, \: \mathcal{S}^k \leftarrow \emptyset, \: \balpha^{(k)} \leftarrow \frac{1}{\mathcal{K}} \mathbf{1}$\;
 $i, \: \by^{(k)} \leftarrow \text{Alg7}(\text{FOP}, \bC, \balpha^{(k)}, \hat\by)$\;
 $\xi_\text{ABS} \leftarrow \bm{\alpha}^{(k) \intercal} C \hat\by - \bm{\alpha}^{(k) \intercal} C \by^{(k)}$\;
 \While{True}{
  $k \leftarrow k + 1$\;
  $\mathcal{S}^k \leftarrow \mathcal{S}^{k-1} \cup \by^{(k-1)}$\;
  $\bm{\alpha}^{(k)} \leftarrow \bm{\alpha}^{(k-1)} - \frac{2}{2+k}(\bm{e}_i - \balpha^{(k-1)})$\;
  $i, \: \by^{(k)} \leftarrow \text{Alg7}(\text{FOP}, \bC, \balpha^{(k)}, \hat\by)$\;
  \eIf{$\balpha^{(k)} = \balpha^{(k-1)}$\textbf{or} $\bm{\alpha}^{(k) \intercal} C \hat{\by} - \bm{\alpha}^{(k) \intercal} C \by^{(k)} > \xi_\text{ABS} + (\balpha^{(k)} - \balpha^{(k-1)})^\intercal (C (\hat{\by} - \by^{(k-1)}))$ \textbf{and} $\by^{(k)} \in \mathcal{S}^k$}{
    $\bm{\alpha}^{\text{final}}, \xi^\text{final} \leftarrow \text{GIO}_{\text{ABS}}(\hat\by, \mathcal{S}^k)$\;
    $k \leftarrow k + 1$\;
    $\by^{(k)} \leftarrow \text{FOP}(\bm{\alpha}^{\text{final}})$\;
    \eIf{$\balpha^{\text{final} \intercal} C \hat\by - \balpha^{\text{final} \intercal} C \by^{(k)} = \xi^\text{final}$}{
    \textbf{stop}
    }{
    $\balpha^{(k)} \leftarrow \balpha^{\text{final}}$\;
    $\mathcal{S}^k \leftarrow \mathcal{S}^{k-1} \cup \by^{(k-1)}$\;
    $t \leftarrow \frac{t}{2}$\;
    }
  }{
   $\xi_\text{ABS} \leftarrow \bm{\alpha}^{(k) \intercal} C \hat{\by} - \bm{\alpha}^{(k) \intercal} C \by^{(k)}$\;
 } 
}
\caption{Frank-Wolfe Partial FOP Inverse Method}
\label{ggfwpartial}
\end{algorithm}

\begin{algorithm}
 \SetKwInOut{Input}{input}\SetKwInOut{Output}{output}
 \Input{$G = \{V, E\}, \: n, \: \eta$}
 \Output{$\mathbb{G}$}
 $\mathbb{G} \leftarrow \{\}$\;
 $\bw_e \leftarrow 1, \:\:\:\:\:\:\: \forall e \in E$\;
 \For{$i \in \text{range}(n)$}{
   $\textbf{o}_e \leftarrow \bw_e \text{Exp}(\lambda = 1), \:\:\:\:\:\:\: \forall e \in E$\;
   $O \leftarrow [e \in E | \textbf{o} \text{ sorted in ascending order}]$\;
   $\mathcal{C} \leftarrow \{\}$ \;
   \For{$e = (v_1, v_2) \in O$}{
     \eIf{$\exists c \in \mathcal{C} | v_1 \in c \text{ or } v_2 \in c$}{
       $\bw_e \leftarrow \eta \bw_e$
     }{
       $\mathcal{C} \leftarrow \mathcal{C} \cup e$\;
       $\bw_e \leftarrow \frac{\bw_e}{\eta}$
     }
   }
   $G' \leftarrow G$\;
   contract all edges $e$ of $G'$ where $e \in \mathcal{C}$\;
   $\mathbb{G} \leftarrow \mathbb{G} \cup G'$
 }
 \caption{Boosted Coarsening Ensemble Generation Method}
 \label{alg:boosting}
\end{algorithm} 

\begin{algorithm}
 \SetKwInOut{Input}{input}\SetKwInOut{Output}{output}
 \Input{$C, \: \hat{\by}$, FOP, $t$, $n_k$, $K$}
 \Output{$\bm{\alpha}^\text{best}$}
 $k = 0, \: \mathcal{S}^k \leftarrow \emptyset, \: \bm{\alpha}^{(k)} \leftarrow \frac{1}{\mathcal{K}} \mathbf{1}$\;
 \For{$i \in 1 \dots n_k$}{
   $\text{FOP}_{i} \leftarrow \text{Coarsen}(\text{FOP})$\;
   $\by^{(k)}_i \leftarrow \text{FOP}_i(\bm{\alpha}^{(k)})$\;
 }
 $\by^{(k)} \leftarrow \underset{\by^{(k)}_i}{\argmin} \balpha^{(k)\intercal}\bC\by^{(k)}_i$\;
 $\xi_\text{ABS} \leftarrow \bm{\alpha}^{(k) \intercal} C \hat{\by} - \bm{\alpha}^{(k) \intercal} C \by^{(k)}$\;
 \While{$k < K$}{
  $k \leftarrow k + 1$\;
  $\mathcal{S}^k \leftarrow \mathcal{S}^{k-1} \cup \by^{(k-1)}$\;
  $\bm{\alpha}^{(k)} \leftarrow \bm{\alpha}^{(k-1)} - t (C (\hat{\by} - \by^{(k-1)}))$\;
  $\bm{\alpha}^{(k)} \leftarrow \text{proj}_{\Delta^\mathcal{K}}(\bm{\alpha}^{(k)})$\;
  \For{$i \in 1 \dots n_k$}{
    $\text{FOP}_{i} \leftarrow \text{Coarsen}(\text{FOP})$\;
    $\by^{(k)}_i \leftarrow \text{FOP}_i(\bm{\alpha}^{(k)})$\;
  }
  $\by^{(k)} \leftarrow \underset{\by^{(k)}_i}{\argmin} \balpha^{(k)\intercal}\bC\by^{(k)}_i$\;
  $\xi_\text{ABS} \leftarrow \bm{\alpha}^{(k) \intercal} C \hat{\by} - \bm{\alpha}^{(k) \intercal} C \by^{(k)}$\;
}
$\balpha^\text{best} \leftarrow \frac{\underset{i \in 0 \dots K}{\sum} (i+1)^2 \balpha^{(i)}}{\underset{i \in 0 \dots K}{\sum} (i+1)^2}$\;
 \textbf{stop}\;
 \caption{Projected gradient descent method with stochastic subgradient estimation.}
 \label{stochasticggproj}
\end{algorithm}

\end{document}